\documentclass[11pt,a4paper,reqno]{article}

\usepackage[utf8]{inputenc}
\usepackage[affil-it]{authblk}
\usepackage{amsfonts, amsmath, amssymb, amsthm, bbm, color, enumerate, graphicx, mathtools, tikz, hyperref, relsize}
\usepackage[numbers,sort]{natbib}
\usepackage[sc]{mathpazo}
\usepackage[T1]{fontenc}
\usepackage{eulervm}
\usepackage[sc]{mathpazo}
\usepackage[a4paper,includeheadfoot, total={6.5in,10in}]{geometry}
\usepackage[title]{appendix}

\usepackage{subcaption} 

\usepackage{todonotes}
\usepackage{tikz-3dplot}
\usetikzlibrary{decorations.text, decorations.markings, mindmap,trees,shadows, decorations.pathreplacing}
\tikzstyle{doublearr}=[latex-latex,red, line width=0.5pt]
\tikzstyle{doublearr2}=[latex-latex,green!80!black, line width=0.5pt]
\usepackage{xcolor}

\newcommand{\QQ}{\mathbf{Q}}
\newcommand{\bQ}{\bar{Q}}

\newcommand{\cE}{\mathcal{E}}

\newcommand{\R}{\mathbbm{R}}

\newcommand{\dif}{\ensuremath{\mbox{d}}}  
\newcommand{\e}{\mathrm{e}}

\newcommand{\dto}{\ensuremath{\xrightarrow{\mathcal{L}}}}  
\newcommand\disteq{\stackrel{\text{\scriptsize d}}{=\joinrel=}}  

\newcommand\E[1]{\mathbbm{E}\left(#1\right)}  
\newcommand{\expt}{\mathbb{E}}
\newcommand\ind[1]{\ensuremath{\mathbbm{1}_{\left[#1\right]}}} 
\newcommand\Pro[1]{\mathbbm{P}\left(#1\right)}  
\newcommand{\prob}{\mathbb{P}}

\newtheorem{theorem}{Theorem}

\newtheorem{lemma}[theorem]{Lemma}
\newtheorem{proposition}[theorem]{Proposition}
\newtheorem{remark}{Remark}

\let\plainqed\qedsymbol
\newcommand{\claimqed}{$\lrcorner$}

\hypersetup{
 colorlinks,
 linkcolor=red,          
 citecolor=blue,       
 filecolor=red,   
 urlcolor=black, 
 pdftitle={},
 pdfauthor={},
 pdfcreator={Debankur Mukherjee, Sayan Banerjee},
 pdfsubject={},
 pdfkeywords={}
}

\usepackage{fancyhdr}
 
\tikzstyle{mybox} = [draw=red, fill=yellow!20, thick, minimum height=.4cm,
    rectangle, rounded corners]
\tikzstyle{fancytitle} =[fill=blue, text=white]
\usetikzlibrary{mindmap,trees,shadows}

\numberwithin{equation}{section}
\numberwithin{theorem}{section}

\tikzstyle{mybox} = [draw=black, thick, minimum height=.6cm,
    rectangle,text centered]
\tikzstyle{fancytitle} =[fill=blue, text=white]

\begin{document}

\title{Join-the-Shortest Queue Diffusion Limit in Halfin-Whitt Regime: Tail Asymptotics and Scaling of Extrema}

\author[1]{Sayan Banerjee\footnote{\texttt{sayan@email.unc.edu}\hspace{1cm}$^\dagger$\texttt{debankur\_mukherjee@brown.edu}}}
\author[2]{Debankur Mukherjee$^\dagger$}
\affil[1]{
University of North Carolina, Chapel Hill}
\affil[2]{
Brown University}

\renewcommand\Authands{, }

\date{\today}

\maketitle 

\begin{abstract}
Consider a system of $N$~parallel single-server queues with unit-exponential service time distribution and a single dispatcher where tasks arrive as a Poisson process of rate $\lambda(N)$. When a task arrives, the dispatcher assigns it to one of the servers according to the Join-the-Shortest Queue (JSQ) policy. Eschenfeldt and Gamarnik (2015) established that in the Halfin-Whitt regime where $(N - \lambda(N)) / \sqrt{N} \to \beta > 0$ as $N \to \infty$, appropriately scaled occupancy measure of the system under the JSQ policy converges weakly on any finite time interval to a certain diffusion process as $N\to\infty$. Recently, it was further established by Braverman (2018) that the convergence result extends to the steady state as well, i.e., stationary occupancy measure of the system converges weakly to the steady state of the diffusion process as $N\to\infty$, proving the interchange of limits result.

In this paper we perform a detailed analysis of the steady state of the above diffusion process. Specifically, we establish precise tail-asymptotics of the stationary distribution and scaling of extrema of the process on large time interval. Our results imply that the asymptotic steady-state scaled number of servers with queue length two or larger exhibits an Exponential tail, whereas that for the number of idle servers turns out to be Gaussian. From the methodological point of view, the diffusion process under consideration goes beyond the state-of-the-art techniques in the study of the steady state of diffusion processes. Lack of any closed form expression for the steady state and intricate interdependency of the process dynamics on its local times make the analysis significantly challenging. We develop a technique involving the theory of regenerative processes that provides a tractable form for the stationary measure, and in conjunction with several sharp hitting time estimates, acts as a key vehicle in establishing the results. The technique and the intermediate results might be of independent interest, and can possibly be used in understanding the bulk behavior of the process.\\\\



\end{abstract}

\section{Introduction}
For any $\beta>0$, consider the following diffusion process
\begin{equation}\label{eq:diffusionjsq}
\begin{split}
Q_1(t) &= Q_1(0) + \sqrt{2} W(t) - \beta t +
\int_0^t (- Q_1(s) + Q_2(s)) \dif s - L(t), \\
Q_2(t) &= Q_2(0) + L(t) - \int_0^t Q_2(s)\dif s
\end{split}
\end{equation}
for $t \geq 0$, where $W$ is the standard Brownian motion, $L$ is
the unique nondecreasing nonnegative process in~$D_\R[0,\infty)$ satisfying
$\int_0^\infty \mathbbm{1}_{[Q_1(t) < 0]} \dif L(t) = 0$,
and $(Q_1(0), Q_2(0)) \in (-\infty, 0] \times [0, \infty)$.
In this paper we establish tail asymptotics of the stationary distribution of the above diffusion process and identify the scaling behavior of $\inf_{0\leq s\leq t}Q_1(s)$ and $\sup_{0\leq s\leq t}Q_2(s)$ for large $t$.
The diffusion process in~\eqref{eq:diffusionjsq} arises as the weak limit of the sequence of scaled occupancy measure of systems under the Join-the-Shortest Queue (JSQ) policy, as the system size (number of servers in the system) becomes large. 
Specifically, consider a system with $N$~parallel identical single-server queues and a single dispatcher.
Tasks with unit-mean exponential service requirements arrive at the
dispatcher as a Poisson process of rate $\lambda(N)$,
and are instantaneously forwarded to one of the servers with the shortest queue length (ties are broken arbitrarily).
For $t\geq 0$, let
$$\QQ^N(t): =
\left(Q_1^N(t), Q_2^N(t), \dots\right)$$ 
denote the
system occupancy measure, where $Q_i^N(t)$ is the number of servers
under the JSQ policy with a queue length of~$i$ or larger,
at time~$t$, including the possible task in service, $i = 1, 2,\dots$.
Now consider an asymptotic regime where the number of servers grows large, and additionally assume that
$$\frac{N - \lambda(N)}{ \sqrt{N}} \to \beta\quad\text{as}\quad N \to \infty$$
for some positive coefficient $\beta > 0$, i.e., the load per server $\lambda(N)/N$ approaches unity as $1 - \beta / \sqrt{N}$, with $\beta > 0$ some positive coefficient. 
In terms of the aggregate traffic load and total service capacity, this scaling corresponds to the so-called Halfin-Whitt heavy-traffic regime which was introduced in the seminal paper~\cite{HW81} and has been extensively studied since. 
The set-up in~\cite{HW81}, as well as the numerous model extensions in the literature (see~\cite{GG13a, GG13b, HW81, LK11, LK12, FKL14, LMZ17}, and the references therein), predominantly concerned a setting with a single centralized queue and server pool (M/M/N), rather than a scenario with parallel queues.
Eschenfeldt and Gamarnik~\cite{EG15} initiated the study of the scaling behavior for parallel-server systems in the Halfin-Whitt heavy-traffic regime.
Define the centered and scaled system occupancy states as $\bar{\QQ}^N(t) =
\big(\bar{Q}_1^N(t), \bar{Q}_2^N(t), \dots\big)$, with 
$$\bar{Q}_1^N(t) = - \frac{N-Q_1^N(t)}{ \sqrt{N}},\qquad \bar{Q}_i^N(t) =\frac{ Q_i^N(t)}{\sqrt{N}},\quad i = 2, 3\dots.$$
The reason why $Q_1^N(t)$ is centered around~$N$ while $Q_i^N(t)$,
$i = 2, \dots$, are not, is because the fraction of servers at time $t$ with a queue
length of exactly one tends to $1$, whereas the fraction of servers with a queue length of two or more tends to zero as $N\to\infty$. 
For each fixed~$N$, $\bar{\QQ}^N$ is a positive recurrent continuous time Markov chain, and there exists a stationary distribution for $\bar{\QQ}^N(t)$ as $t \rightarrow \infty$. 
Denote by $\bar{\QQ}^N(\infty)$ a random variable distributed as the steady state of the process $\bar{\QQ}^N(t)$.
Assuming $(\bQ_i^N(0))_{i\geq 1} \dto (Q_i(0)))_{i\geq 1}$ with  $Q_i(0)= 0$ for $i \ge 3$, it was shown by Eschenfeldt and Gamarnik~\cite{EG15} that on any finite time interval $[0,T]$,
the sequence of processes $\big\{(\bar{Q}_1^N(t), \bar{Q}_2^N(t), \ldots)\big\}_{0 \le t \le T}$ 
converges weakly to the limit
 $\big\{(Q_1(t), Q_2(t),\ldots)\big\}_{0 \le t \le T}$, where $(Q_1, Q_2)$ is given by~\eqref{eq:diffusionjsq} and $Q_i(\cdot)\equiv 0$ for~$i\geq 3$. 
 Subsequently, a broad class of other schemes were shown to exhibit the same scaling behavior in this regime~\cite{MBLW16-3, MBLW15, MBL17}.
 See~\cite{BBLM18} for a recent survey.
In all these above works, the convergence of scaled occupancy measure was established in the transient regime on any finite time interval.
The tightness of diffusion-scaled occupancy measure and the interchange of limits were open until recently, when Braverman~\cite{Braverman18} further established that the weak convergence result extends to the steady state as well, i.e., $\bar{\QQ}^N(\infty)$ converges weakly to $(Q_1(\infty), Q_2(\infty), 0, 0,\ldots)$ as $N\to\infty$, where $(Q_1(\infty), Q_2(\infty))$ is distributed as the stationary distribution of the process $(Q_1, Q_2)$.
 Thus, the steady state of the diffusion process in~\eqref{eq:diffusionjsq} captures the asymptotic behaviors of large-scale systems under the JSQ policy.
 
The steady-state of the diffusion process in~\eqref{eq:diffusionjsq} is technically hard to analyze.
In fact, even establishing its ergodicity is non-trivial. The standard method employed in studying steady-state behavior of diffusions \cite{ABD, BL1,DW, HM} is to construct a suitable Lyapunov function which shows that the diffusion has a strong drift towards a compact set. Inside the compact set, some irreducibility condition, like uniform ellipticity (as in \cite{ABD,BL1,DW}) or hypoellipticity (as in \cite{HM}), is used to show positive recurrence, and consequently, existence and uniqueness of the stationary distribution and ergodicity of the diffusion process. The construction of the Lyapunov function usually involves establishing stability of the associated noiseless dynamical system and having tractable bounds on hitting times for this deterministic system. In our setup, even the noiseless system requires non-trivial analysis (see Section 4.1 of \cite{Braverman18}).
In~\cite{Braverman18} a Lyapunov function is obtained via a generator expansion framework using the Stein's method that establishes exponential ergodicity of $(Q_1,Q_2)$.
Although this approach gives a good handle on the rate of convergence to stationarity, the non-trivial dynamics of the noiseless system results in a complicated form for the Lyapunov function which sheds little light on the form of the stationary distribution itself. 
Moreover, the diffusion in~\eqref{eq:diffusionjsq} (without the reflection term) is not hypoelliptic and this complicates things even further. It is also worth pointing out here that we obtain different tail behavior for $Q_1$ and $Q_2$ (Gaussian and Exponential, respectively) and get explicit dependence of $\beta$ in the exponents, which is hard to obtain using the Lyapunov function methods known in the literature.

This demands a fundamentally different characterization of the stationary distribution.
For that we take resort to the theory of regenerative processes (see Chapter 10 of \cite{Thorisson}) to obtain a tractable representation of the steady state.
A variant of this method was first used in~\cite{banerjee2015} to study a diffusion process with inert drift, although the stationary distribution in that case had an explicit product form that facilitated the analysis, as opposed to the current scenario.
First, we show that the diffusion $\{(Q_1(t), Q_2(t))\}_{t \ge 0}$ can be decomposed into i.i.d.~renewal cycles between carefully constructed regeneration times having good moment bounds. 
This decomposition gives an alternative, more transparent proof of ergodicity, and also shows that the diffusion falls in the category of classical regenerative processes. 
Loosely speaking, regeneration times are random times when the process \emph{starts afresh}, and the theory of classical regenerative processes can be used to conclude that the stationary behavior of a process is same as the behavior within one renewal cycle (i.e., between two successive regeneration times).
The regenerative process representation enables us to obtain a form for the stationary distribution that is amenable to analysis (see Theorem \ref{th:stationary}).
Tail estimates for the stationary measure are then obtained by analyzing this form and are presented in Theorem \ref{th:statail}. 
Moreover, in Theorem~\ref{th:lil}, we obtain precise almost sure
scaling behavior of the extrema of the process sample paths.

The regenerative structure of the diffusion process and the intermediate results might be of independent interest.
In fact, they might also be used to provide detailed result about the behavior of the stationary measure near the center (bulk behavior) and produce sharp estimates on the stationary mean of $Q_2$. \\
 
Rest of the article is arranged as follows.
 In Section~\ref{sec:main}, we describe the two main results of this paper. 
 In the Section~\ref{sec:reg}, we establish $(Q_1, Q_2)$ as a classical regenerative process and state several crucial hitting time estimates that are required to prove the main results.
 In Section~\ref{sec:renewal}, we obtain a tail estimate for the regeneration time which, in particular, implies that it has a finite first moment. This, in turn, implies the ergodicity of the diffusion process and gives a tractable form for the stationary distribution. 
 In Section~\ref{sec:fluc}, we obtain fluctuation estimates of the paths of $Q_1$ and $Q_2$ between two successive regeneration times, which are used in the proofs of Theorems~\ref{th:statail} and~\ref{th:lil}. 
 In Section~\ref{sec:proofmain}, we combine the results in Sections~\ref{sec:reg}, \ref{sec:renewal} and~\ref{sec:fluc} to prove Theorems~\ref{th:statail} and~\ref{th:lil}.

\section{Main results}\label{sec:main}

In this section we will state the main results, and discuss their ramifications. 
Recall the diffusion process $\{(Q_1(t), Q_2(t))\}_{t \ge 0}$ as defined by Equation~\eqref{eq:diffusionjsq}.
As mentioned in the introduction, it is known~\cite{Braverman18} that for any $\beta>0$, $(Q_1, Q_2)$ is an ergodic continuous-time Markov process.
Let $(Q_1(\infty), Q_2(\infty))$ denote a random variable distributed as the unique stationary distribution~$\pi$ of the process.
Then the next theorem gives a precise characterization of the tail of the stationary distribution.
\begin{theorem}\label{th:statail}
For any $\beta>0$ there exist positive constants $C_1, C_2, D_1, D_2$ \emph{not} depending on $\beta$ and positive constants $C^{l}(\beta), C^{u}(\beta), D^{l}(\beta), D^{u}(\beta), C_R(\beta), D_R(\beta)$ depending \emph{only} on $\beta$ such that 
\begin{equation}
\begin{split}
C^{l}(\beta)\e^{-C_1x^2} \le \pi(Q_1(\infty) < -x) \le C^{u}(\beta)\e^{-C_2x^2}, \ \ x \ge C_R(\beta)\\
D^{l}(\beta)\e^{-D_1\beta y} \le \pi(Q_2(\infty) > y) \le D^{u}(\beta)\e^{-D_2\beta y}, \ \ y \ge D_R(\beta).
\end{split}
\end{equation}
\end{theorem}
The dependence on $\beta$ of the tail-exponents is precisely captured in the above theorem.
Note that $Q_1(\infty)$ has a Gaussian tail, and the tail exponent is uniformly bounded by constants which do not depend on $\beta$, whereas $Q_2(\infty)$ has an exponentially decaying tail, and the coefficient in the exponent is linear in $\beta$.
Loosely speaking, Theorem \ref{th:statail} implies that the sample path of $Q_2$ tends to spend more time taking larger values as $\beta$ becomes smaller, whereas the sample path of $Q_1$ seems to be less affected by $\beta$.
Also, note that the dependence of the exponents on $\beta$ is useful in obtaining the growth rate of the extreme values of $Q_1$ and $Q_2$ on large time intervals, as further made precise in Theorem~\ref{th:lil} below.
\begin{remark}\label{rem:comp1}\normalfont
Let us now discuss a further implication of Theorem~\ref{th:statail}.
Recall that $Q_i^N(t)$ denotes the number of servers in the $N$-th system with queue length $i$ or larger at time $t$.
Let $S^N(t):=\sum_{i\geq 1}Q_i^N(t)$ denote the total number of tasks in the system.
Then \cite[Theorem 5]{Braverman18} implies that $(S^N(\infty)-N)/\sqrt{N}$ converges weakly to $S(\infty)\disteq Q_1(\infty)+Q_2(\infty)$.
In that case, Theorem~\ref{th:statail} implies that $S(\infty)$ has an Exponential upper tail (large positive deviation) and a Gaussian lower tail (large negative deviation).
\end{remark}

\begin{remark}\label{rem:comp2}\normalfont
It is worth mentioning that in case of M/M/N systems in the Halfin-Whitt heavy-traffic regime~\cite[Theorem 2]{HW81}, the centered and scaled total number of tasks in the system $(\bar{S}^N(t)-N)/\sqrt{N}$ converges weakly to a diffusion process $\{\bar{S}(t)\}_{t\geq 0}$ having the infinitesimal generator 
$A = (\sigma^2(x)/2)(\dif^2/\dif x^2)+m(x)(\dif/\dif x)$
with
\[m(x) = \begin{cases}
-\beta & \mbox{ if } x >0\\
-(x+\beta) & \mbox{ if } x \le 0
\end{cases}\qquad \text{and}\qquad\sigma^2(x) = 2.\]
Note that since this is a simple combination of a Brownian motion with a negative drift (when all servers are fully occupied) and an Ornstein Uhlenbeck process (when there are idle servers), the steady-state distribution $\bar{S}(\infty)$ can be {\em computed explicitly}, and is a combination of Exponential (from the Brownian motion with a negative drift) and Gaussian (from the OU process).
Although in terms of tail asymptotics, $\bar{S}(\infty)$ and $S(\infty)$ in Remark~\ref{rem:comp1} behave somewhat similarly, there are some fundamental differences between the two processes, that not only make the analysis of the JSQ policy much harder, but also lead to several completely different qualitative behavior.
\begin{enumerate}[{\normalfont (i)}]
\item Observe that in case of M/M/N, whenever there are some waiting tasks (equivalent to $Q_2$ being positive in our case), the queue length has a constant negative drift towards zero. This leads to the Exponential upper tail of $\bar{S}(\infty)$, by comparing with the stationary distribution of reflected Brownian motion with constant negative drift. In our case, the rate of decrease of $Q_2$ is always proportional to itself, which makes it somewhat counter-intuitive that its stationary distribution has an Exponential tail.
\item Further, from \eqref{eq:diffusionjsq}, $Q_2$ {\em never} hits zero. 
Thus, in the steady state, there is no mass at $Q_2=0$, and the system always has waiting tasks. 
This is in sharp contrast with the M/M/N case, where with positive probability the steady-state system has no waiting task. 
\item In the M/M/N setup, given that a task faces a non-zero wait, \emph{the steady-state waiting time is of order $1/\sqrt{N}$} unlike in our case, where it is of \emph{constant order} (the time till the service of the task ahead of it in its queue finishes). 
Moreover, in the current scenario, it is easy to see that $Q_1$ (the limit of the scaled number of idle servers) spends zero time at the origin, i.e., in steady state  the fraction of arriving tasks that find all servers busy vanishes in the large-N limit.
Consequently, JSQ achieves an {\em asymptotically vanishing steady-state probability of non-zero wait} (in fact, this is of order $1/\sqrt{N}$, see~\cite{Braverman18}). 
This is another sharp contrast with the M/M/N case, where the \emph{asymptotic steady-state probability of non-zero wait is strictly positive}.
\item In the M/M/N setup, the number of idle servers can be non-zero only when the number of waiting tasks is zero.
Thus, the dynamics of both the number of idle servers and the number of waiting tasks are completely captured by the one-dimensional process $S^N$ and by the one-dimensional diffusion  $\bar{S}$ in the limit. 
But in our case, $Q_2$ is never zero, and the dynamics of $(Q_1, Q_2)$ is truly two-dimensional (although the diffusion is non-elliptic) with $Q_1$ and $Q_2$ interacting with each other in an intricate manner.
\end{enumerate}
\end{remark}

The next theorem establishes scaling behavior of the extrema of the process $\{(Q_1(t), Q_2(t))\}_{t \ge 0}$ on large time intervals.
\begin{theorem}\label{th:lil}
There exists a positive constant $\mathcal{C^*}$ not depending on $\beta$ such that the following hold almost surely along any sample path:
\begin{align*}
-2\sqrt{2} &\le \liminf_{t \rightarrow \infty} \frac{Q_1(t)}{\sqrt{\log t}} \le -1,\\
\frac{1}{\beta} &\le \limsup_{t \rightarrow \infty} \frac{Q_2(t)}{\log t} \le \frac{2}{\mathcal{C^*} \beta}.
\end{align*}
\end{theorem}
Again, Theorem~\ref{th:lil} captures the explicit dependence on $\beta$ of the width of the fluctuation window of $Q_1$ and $Q_2$. 
Specifically, note that the width of fluctuation of $Q_1$ does not depend on the value of $\beta$, whereas that of $Q_2$ 
is linear in $\beta^{-1}$. 
\begin{remark}\normalfont
From the proofs of Theorems~\ref{th:statail} and~\ref{th:lil} one can see that $C_1=4$, $C_2=1/16$, $D_1=2$, $D_2=1/16$, and $\mathcal{C^*} = 1/16$.
However, we are not explicit about them in the statements of the theorems since these estimates are not sharp in the constants.
\end{remark}
\section{Regenerative process view of the diffusion}\label{sec:reg}
As mentioned in the introduction, the key challenge in analyzing the steady state of the diffusion process in~\eqref{eq:diffusionjsq} stems from its lack of explicit characterization. 
In order to obtain sharp estimates for the stationary distribution we take resort to the theory of regenerative processes.
Loosely speaking, a stochastic process is called \emph{classical regenerative} if it starts anew at random times (called \emph{regeneration times}), independent of the past. See~\cite[Chapter 10]{Thorisson} for a rigorous treatment of regenerative processes.
The regeneration times split the process into renewal cycles that are independent and identically distributed, possibly except the first cycle.
Consequently, the behavior inside a specific renewal cycle characterizes the steady-state behavior.

In case of recurrent discrete state-space Markov chains regeneration times can be defined as hitting times of a fixed state.
Although the diffusion process in~\eqref{eq:diffusionjsq} is two dimensional, we will show that it actually \textit{exhibits point recurrence} and we can define regeneration times in terms of hitting times as follows.

First we introduce the following notations.
\begin{align*}
\tau_i(z)&:=\inf\{t \ge 0: Q_i(t)=z\}, \ \ i=1,2.\qquad\text{and}\qquad
\sigma(t):= \inf\{s \ge t: Q_1(s)=0\}.
\end{align*}
We now define the renewal cycles as follows.
Fix any $B>0$. For $k\geq 0$, define the stopping times
\begin{align}\label{rendef}
\alpha_{2k+1} &:= \inf \Big\{t\geq \alpha_{2k} : Q_2(t) = B\Big\},\quad
\alpha_{2k+2} := \inf \left\{t>\alpha_{2k+1}:Q_2(t) = 2B\right\},\quad
\Xi_k := \alpha_{2k+2},
\end{align}
with the convention that $\alpha_0=0$ and $\Xi_{-1}=0$.
The dependence of $B$ in the above stopping times is suppressed for convenience in notation. 
Hereafter we will assume $B>0$ to be fixed unless mentioned otherwise.
The next lemma describes the diffusion process as an appropriate classical regenerative process.
\begin{lemma}\label{lem:regeneration}
The process $\{Q_1(t), Q_2(t)\}_{t\geq 0}$ is a classical regenerative process with regeneration times given by $\{\Xi_k\}_{k\geq 0}$.
\end{lemma}
\begin{proof}
Note that it is enough to prove that 
$Q_1(\alpha_{2k}) = 0$ for all $k\geq 1$.
Indeed, this ensures that for all $k\geq 0$, $(Q_1(\Xi_k), Q_2(\Xi_k)) = (0, 2B)$, and the Markov process naturally regenerates at time $\Xi_k$.

Fix any $k\geq 1$.
Assume, if possible, $Q_1(\alpha_{2k}) < 0$. 
In that case, the path-continuity of $Q_1$ implies that the local time $L$ is constant in a small neighborhood of $\alpha_{2k}$.
Consequently, $Q_2$ must be strictly decreasing in an open time interval containing $\alpha_{2k}$.
This contradicts the fact that $\alpha_{2k}$ is the hitting time of a level \emph{from below} by the process $Q_2$.
\end{proof}
The above lemma implies that the regenerative cycles given by $\{(Q_1(t), Q_2(t))\}_{\Xi_k\leq t<\Xi_{k+1}}$ form an i.i.d.~sequence for $k\geq 0$.
The time intervals $\{\Xi_{k+1}-\Xi_k\}_{k\geq 0}$ are called the \emph{inter-regeneration times}.
In order to characterize the steady-state distribution using regenerative approach, we first show that the initial \emph{delay length} $\Xi_0$ (time to enter into the regenerative cycles starting from an arbitrary state) as well as inter-regeneration times have finite expectations.
In fact, the next proposition establishes detailed tail asymptotics for the delay length $\Xi_0$ and thus, in particular, for the inter-regeneration times.
\begin{proposition}\label{prop:Xi}
Let $(Q_1(0), Q_2(0))= (x,y)$ with $x \le 0, y >0$. 
There exist constants $ c^{(1)}_\Xi, c^{(2)}_\Xi, t_{\Xi}>0$, possibly depending on $x,y,B,\beta$, such that for all $t\geq t_{\Xi}$,
$$\prob_{(x,y)}(\Xi_0>t)\leq c^{(1)}_\Xi\exp(-c^{(2)}_\Xi t^{1/6}).$$
In particular, $\expt_{(x,y)}\Xi_0<\infty.$
\end{proposition}
Proposition~\ref{prop:Xi} is proved in Section~\ref{sec:renewal}.
Proposition \ref{prop:Xi} yields the existence and uniqueness of the stationary distribution and ergodicity of the process as stated in Theorem~\ref{th:stationary} below. 
We note that the geometric ergodicity has already been proved in~\cite{Braverman18}.
The principal importance of Theorem~\ref{th:stationary} lies in the fact that it provides an explicit form of the stationary measure which will be the key vehicle in the study of the tail asymptotics and the fluctuation window, as stated in Theorems~\ref{th:statail} and~\ref{th:lil}. 
\begin{theorem}\label{th:stationary}
Fix any $B>0$. The process described by Equation~\eqref{eq:diffusionjsq} has a unique stationary distribution $\pi$ which can be represented as
$$
\pi((Q_1(\infty), Q_2(\infty)) \in A) = \dfrac{\mathbb{E}_{(0, 2B)}\left(\int_{0}^{\Xi_0}\mathbbm{1}_{[(Q_1(s),Q_2(s)) \in A]}ds\right)}{\mathbb{E}_{(0, 2B)}\left(\Xi_0\right)}
$$
for any measurable set $A \subseteq (-\infty, 0] \times (0, \infty)$. Moreover, the process is ergodic in the sense that for any measurable function $f$ satisfying $\mathbb{E}_{(0, 2B)}\left(\int_{0}^{\Xi_0}f((Q_1(s),Q_2(s)))ds\right) < \infty$,
\begin{equation}\label{ergodicity}
\frac{1}{t}\int\limits_0^t f((Q_1(s),Q_2(s)))ds \longrightarrow  \frac{\mathbb{E}_{(0, 2B)}\left(\int_{0}^{\Xi_0}f((Q_1(s),Q_2(s))ds\right)}{\mathbb{E}_{(0, 2B)}\left(\Xi_0\right)}
\end{equation}
almost surely as $t \rightarrow \infty$.
\end{theorem}
The above theorem follows using \cite[Chapter 10, Theorem 2.1]{Thorisson}, details of which are deferred till Section~\ref{sec:renewal}.
\begin{remark}\normalfont
We note that it can be shown by soft arguments involving Girsanov theorem and the theory of L\'evy processes that the distribution of $\Xi_1 - \Xi_0$ has a density with respect to the Lebesgue measure, see the proof of Lemma 7.1 in \cite{banerjee2015}. 
This implies that the inter-regeneration time $\Xi_{k+1}-\Xi_k$ is \emph{spread-out} (see Section 3.5 of Chapter 10 in~\cite{Thorisson}).
Consequently, the total variation convergence of the diffusion process at time $t$ to the stationary distribution as $t\to\infty$, can be obtained using Theorem 3.3 of Chapter 10 in~\cite{Thorisson}. 
However, we skip this argument, since geometric ergodicity has already been established in \cite[Theorem 3]{Braverman18}. 
\end{remark}
In light of Theorem~\ref{th:stationary}, observe that establishing tail asymptotics of the stationary distribution reduces to studying the amount of time spent by the diffusion in a certain region in one particular renewal cycle.
The next theorem provides several important hitting time estimates that will play a crucial role in the proofs of Theorems~\ref{th:statail} and~\ref{th:lil}.
Define
\begin{equation}\label{lzerodef}
l_0(\beta) := \max\left\lbrace \beta, \beta^{-1}, \frac{1}{\beta}\log\frac{1}{\beta}\right\rbrace.
\end{equation}
\begin{theorem}\label{th:excren}
There exists a positive constant $R_0$ such that with $B = R_0l_0(\beta)$ in~\eqref{rendef},
the following hold:
\begin{itemize}
\item[{\normalfont(i)}] There exist constants $C^*_1, C^*_2 >0$ that do not depend on $\beta$ such that for all $y \ge 4B$,
\begin{equation*}
\prob_{(0, 2B)}\left(\tau_2(y) \le \Xi_0 \right) \le C^*_1 \e^{-C^*_2 \beta (y-\beta)/2}.
\end{equation*}
\item[{\normalfont(ii)}] 
For all $y \ge 2B$,
\begin{equation*}
\prob_{(0,2B)}\left(\tau_2(y) \le \Xi_0\right) \ge (1-\e^{-\beta R_0l_0(\beta)})\e^{-\beta(y-2R_0l_0(\beta))}.
\end{equation*}
\item[{\normalfont(iii)}] There exists a constant $C^*(\beta)>0$ depending on $\beta$ such that for any $x \ge 18B$,
\begin{equation*}
\prob_{(0,2B)}\left(\tau_1(-x) \le \Xi_0\right) \le C^*(\beta) \e^{-(x-2\beta)^2/8}.
\end{equation*}
\item[{\normalfont(iv)}] There exists a constant $C^{**}(\beta) >0$ depending on $\beta$ such that for any $x \ge \beta$,
\begin{equation*}
\prob_{(0,2B)}\left(\inf_{t \le \Xi_0} Q_1(t) < -x\right) \ge C^{**}(\beta) \e^{-x^2}.
\end{equation*}
\end{itemize}
\end{theorem}
Theorem~\ref{th:excren} is proved in Section~\ref{sec:fluc} where we analyze the behavior of the process $(Q_1, Q_2)$ between two successive regeneration times. 
Results in Theorem~\ref{th:excren} in conjunction with Proposition~\ref{prop:Xi} and Theorem~\ref{th:stationary} are used to prove Theorems~\ref{th:statail} and~\ref{th:lil}, which is presented in Section~\ref{sec:proofmain}.

\section{Analysis of regeneration times}\label{sec:renewal}
In this section we will prove Proposition~\ref{prop:Xi} and Theorem~\ref{th:stationary}.
The proof of Proposition~\ref{prop:Xi} consists of several steps.
The first step is to analyze the down-crossings of $Q_2$, where 
we establish various hitting time estimates in the time interval $[\alpha_{2k}, \alpha_{2k+1}]$, $k\geq 0$.
In particular, we prove the following lemma.
\begin{lemma}\label{lem:alphadowntail}
Fix $(Q_1(0), Q_2(0))=(x,y)$ with $x \le 0, y>0$. 
There exist $c_{\alpha_1}, c_{\alpha_1}', t_{\alpha_1}>0$ possibly depending on $(x,y)$, $B$, and $\beta$, such that for all $t \ge t_{\alpha_1}$,
$$
\prob_{(x,y)}(\alpha_1 >t) \le c_{\alpha_1}'\exp(-c_{\alpha_1}t^{1/6}).
$$
\end{lemma}
As before, note that setting $(x,y) = (0,2B)$ furnishes the corresponding probabilities when $\alpha_1$ is replaced by $\alpha_{2k+1} - \alpha_{2k}$.
Lemma~\ref{lem:alphadowntail} is proved in Subsection~\ref{ssec:down}.
Next we consider the up-crossings of $Q_2$, where 
we establish various hitting time estimates in the time interval $[\alpha_{2k+1}, \alpha_{2k+2}]$, $k \ge 0$. 
Specifically, we establish the following.
\begin{lemma}\label{lem:alphauptail}
Fix $(Q_1(0), Q_2(0))=(x,y)$ with $x \le 0, y>0$. 
There exist $c_{\alpha_2}, c_{\alpha_2}', t_{\alpha_2}>0$ possibly depending on $(x,y)$, $B$, and $\beta$, such that for all $t \ge t_{\alpha_2}$,
$$
\prob_{(x,y)}(\alpha_2 - \alpha_1 >t) \le c_{\alpha_2}'\exp(-c_{\alpha_2}t^{1/6}).
$$
\end{lemma}
Lemma~\ref{lem:alphauptail} is proved in Subsection~\ref{ssec:up}.
Now observe that Lemmas~\ref{lem:alphadowntail} and~\ref{lem:alphauptail} together complete the proof of Proposition~\ref{prop:Xi}.\hfill $\qed$

\begin{proof}[Proof of Theorem~\ref{th:stationary}]
Due to Proposition~\ref{prop:Xi}, the fact that $\pi$ defined in the theorem is stationary follows from \cite[Chapter 10, Theorem 2.1]{Thorisson}. Now, we will prove the ergodicity result \eqref{ergodicity} which will also yield uniqueness. Take any starting point $(x,y)$ with $x \le 0$ and $y>0$ and recall $\Xi_{-1}=0$. Take any measurable function $f$ satisfying $\mathbb{E}_{(0, 2B)}\left(\int_{0}^{\Xi_0}f((Q_1(s),Q_2(s)))ds\right) < \infty$. Let $N_t = \sup\{k \ge -1: \Xi_{k} \le t\}$. Assume without loss of generality that $f$ is non-negative (for general $f$, consider the positive and negative parts of $f$ separately). We can write
\begin{multline*}
\int_0^{\Xi_0 \wedge t}f((Q_1(s),Q_2(s)))ds + \mathbbm{1}_{[\Xi_1 \le t]}\sum_{k=1}^{N_t}\int_{\Xi_{k-1}}^{\Xi_{k}} f((Q_1(s),Q_2(s)))ds 
\le \int_0^t f((Q_1(s),Q_2(s)))ds\\
\le \int_0^{\Xi_0}f((Q_1(s),Q_2(s)))ds + \sum_{k=1}^{N_t+1}\int_{\Xi_{k-1}}^{\Xi_{k}} f((Q_1(s),Q_2(s)))ds.
\end{multline*}
Clearly, $t^{-1}\int_0^{\Xi_0}f((Q_1(s),Q_2(s)))ds \rightarrow 0$ as $t \rightarrow \infty$. By Proposition 7.3 of \cite{Ross},
$$
t^{-1} \sum_{k=1}^{N_t}\int_{\Xi_k}^{\Xi_{k+1}} f((Q_1(s),Q_2(s)))ds \rightarrow  \frac{\mathbb{E}_{(0, 2B)}\left(\int_{0}^{\Xi_0}f((Q_1(s),Q_2(s))ds\right)}{\mathbb{E}_{(0, 2B)}\left(\Xi_0\right)}
$$
and
 $$
t^{-1} \sum_{k=1}^{N_t+1}\int_{\Xi_k}^{\Xi_{k+1}} f((Q_1(s),Q_2(s)))ds \rightarrow  \frac{\mathbb{E}_{(0, 2B)}\left(\int_{0}^{\Xi_0}f((Q_1(s),Q_2(s))ds\right)}{\mathbb{E}_{(0, 2B)}\left(\Xi_0\right)}
$$
almost surely as $t \rightarrow \infty$. This proves \eqref{ergodicity}, and consequently uniqueness of the stationary distribution.
\end{proof}

\subsection{Down-crossings of \texorpdfstring{$\mathbf{Q_2}$}{Q2} and tightness estimates}\label{ssec:down}
In this subsection, we will prove tail-asymptotics for the distribution of $\alpha_1$ as stated in Lemma~\ref{lem:alphadowntail}. 
This will require a crucial tightness estimate for the process $Q_2$, which is given in Lemma \ref{lem:q2regeneration} below.

\begin{lemma}\label{lem:q2regeneration}
There exist positive constants $c'_1, c'_2, c'_3, c'_4$ not depending on $\beta$ such that the following hold: 
\begin{enumerate}[{\normalfont (i)}]
\item For $\beta \ge 1$ and any $y \ge 1$, for all $ t\geq c'_4 y/\beta$
\begin{align*}
\prob_{(0,\ y + c'_1\beta)}\big(\inf_{s\leq t}Q_2(s) > c'_1\beta \big)\leq c'_3\exp(-c'_2\beta^{2/5}t^{1/5}).
\end{align*}
\item For $\beta \in (0,1)$ and any $y \ge 1$, for all $t\geq c'_4 \big(y\beta^{-1} \vee \beta^{-2}\big)$
\begin{align*}
&\prob_{(0,\ y + c'_1\beta^{-1})}\Big(\inf_{s\leq t}Q_2(s) > \frac{c'_1}{\beta}\Big)
\leq c'_3\left(\exp(-c'_2\beta^{-\frac{2}{5}}t^{\frac{1}{5}}) + \exp(-c'_2\beta^2 t) + \beta^{-2}\exp(-c'_2 t)\right).
\end{align*}
\end{enumerate}
\end{lemma}

Lemma~\ref{lem:q2regeneration} is proved in Appendix~\ref{app:lem4.3}.
In the proof of Lemma~\ref{lem:q2regeneration}, we need to have fine estimates for the time  $Q_2$ takes to hit the level $B$ starting from a large initial state.
This, in turn, amounts to estimating the time integral of the $Q_1$ process when $Q_2$ is large.
The estimate for the time integral, along with several tail probability estimates, completes the proof of Lemma~\ref{lem:q2regeneration}.

We now proceed to prove Lemma~\ref{lem:alphadowntail}.
\begin{proof}[Proof of Lemma~\ref{lem:alphadowntail}]
From Lemma \ref{lem:q2regeneration}, for any $\beta>0$, we obtain $M^*>2B, t^*>0$ such that for all $t \ge t^*$,
\begin{align}\label{down}
\prob_{(0,2M^*)}(\tau_2(M^*) >t) \le C_1\exp(-C_2t^{1/5}),
\end{align}
where the constants $C_1, C_2>0$ depend on $\beta,M^*$.  
Set the starting state to be $(Q_1(0), Q_2(0))=(0,y)$ where $M^* \ge y \ge B$. 
It will be clear from the proof that the same argument works for general starting points $(x,y)$ with $x \le 0, y >0$. For $k \ge 0$, define the following stopping times:
\begin{align*}
\alpha^*_{2k+1} & = \inf \Big\{t\geq \alpha^*_{2k} : Q_2(t) = 2M^* \text{ or } Q_2(t) = B\Big\},\\
\alpha^*_{2k+2} & = \inf \left\{t>\alpha^*_{2k+1}:Q_2(t) = M^* \text{ or } Q_2(t) = B\right\},
\end{align*}
where by convention, we take $\alpha^*_0=0$. 
Let $\mathcal{N}' := \inf\ \{k \ge 0: Q_2(\alpha^*_{2k}) = B\}$.

We will first prove the following: for some positive constant $p(M^*)$ that depends only on $M^*$
\begin{equation}\label{alphadown}
\inf_{z \in [B, M^*]}\prob_{(0,z)}(\tau_2(B) < \tau_2(2M^*))
 \ge p(M^*)>0.
\end{equation}
To see this, recall $S(t) = Q_1(t) + Q_2(t)$ and note that for $t \le \tau_1(-\beta/2)$,
$$
S(0) + \sqrt{2}W(t) - \beta t/2 \ge S(t) \ge Q_1(t).
$$
Further, note that $S(t) \le Q_2(t)$. 
Moreover, due to arguments similar to Lemma~\ref{lem:regeneration}, we know $Q_1(\tau_2(2M^*))=0$, and hence, $S(\tau_2(2M^*)) = Q_2(\tau_2(2M^*))$. 
Combining these facts, we obtain for any $z \in [B, M^*]$,
\begin{equation}\label{eq:lem3.1-p1m}
\begin{split}
\prob_{(0,z)}(\tau_2(2M^*) \le \tau_1(-\beta/2)) &\le \prob(S(t) \text{ hits } 2M^* \text{ before } -\beta/2)\\
&\le \prob(S(0) + \sqrt{2}W(t) - \beta t/2 \text{ hits } 2M^* \text{ before } -\beta/2)\\
&\le \prob(\sqrt{2}W(t) - \beta t/2 \text{ hits } M^* \text{ before } -(M^* + \beta/2))\\
&\le \e^{-\beta M^*/2} <1,
\end{split}
\end{equation}
where we used the fact that the scale function (see \cite[V.46]{RW2000_2}) for $\sqrt{2}W(t) - \beta t/2$ is $s(x)= \exp(\beta x/2)$.\\

Now we will show that if the process $(Q_1, Q_2)$ starts with the initial state $(-\beta/2,z)$ with $z \le 2M^*$, 
then with positive probability $Q_1(t) <0$ for all $t \le \log (2M^* B^{-1})$.  
This in turn implies that $Q_2$ hits the level $B$ before time $\log (2M^* B^{-1})$, since for $t \le \tau_1(0)$, $(\dif/\dif t) Q_2(t) = - Q_2(t)$. 

Construct the Ornstein-Uhlenbeck process $Q_1^+$ on the same probability space as $Q_1$ as follow
$$
Q_1^+(t) = Q_1(0) + \sqrt{2}W(t) + \int_0^t(-Q_1^+(s) + (2M^*-\beta))ds,
$$
where the driving Brownian motion $W$ is the same as that for $Q_1$. 
By~\cite[Proposition 2.18]{Karatzas}, $Q_1(t) \le Q_1^+(t)$ for all $t \le \tau_1(0)$.
Now define the following event
$$
\cE(M^*) := \Big\{Q_1^+(t) < 0 \text{ for all } t \le \log (2M^* B^{-1})\Big\}.
$$
Note that $\cE(M^*)$ does not depend on $z$. It follows from the Doob representation for Ornstein-Uhlenbeck processes that $\prob(\cE(M^*))>0$.
Thus,
\begin{equation}\label{timelength}
\begin{split}
&\inf_{z \in [B, 2M^*)}\prob_{(-\beta/2,z)}\Big(\tau_2(B) \le  \log (2M^* B^{-1}) < \tau_2(2M^*)\Big)\\
&\hspace{2cm}\ge \inf_{z \in [B, 2M^*)}\prob_{(-\beta/2,z)}(\tau_1(0) \ge \log (2M^* B^{-1})) \ge \prob(\mathcal{E}(M^*))>0.
\end{split}
\end{equation}
The strong Markov property in combination with \eqref{eq:lem3.1-p1m} and~\eqref{timelength} now produces the following bound
\begin{align*}
\inf_{z \in [B, M^*]}\prob_{(0,z)}(\tau_2(B) < \tau_2(2M^*))
 \ge (1-\e^{-\beta M^*/2})\prob(\mathcal{E}(M^*))>0
\end{align*}
which proves \eqref{alphadown}. By virtue of \eqref{alphadown}, we have the following for $n \ge 1$,
\begin{align}\label{down1}
\prob(\mathcal{N}'>n) \le (1-p(M^*))^n.
\end{align}
Now, let $T(M^*)$ be a number large enough such that
\begin{equation}\label{downadd}
\prob\left(\sqrt{2} W(T(M^*)) \ge \beta T(M^*)/2 - (2M^* + \beta/2)\right) \le \prob(\mathcal{E}(M^*))/2.
\end{equation}
Then,
\begin{equation}\label{downadd1}
\begin{split}
&\sup_{z \in [B, 2M^*)} \prob_{(0,z)}\left(\tau_2(B) \wedge \tau_2(2M^*) > T(M^*) + \log (2M^* B^{-1})\right)\\
&\le \sup_{z \in [B, 2M^*)} \prob_{(0,z)}\left(\tau_1(-\beta/2) < T(M^*), \tau_2(B) \wedge \tau_2(2M^*) > T(M^*)+ \log (2M^* B^{-1})\right)\\
&\hspace{7cm}+ \sup_{z \in [B, 2M^*)} \prob_{(0,z)}\left(\tau_1(-\beta/2) \ge T(M^*)\right)\\ 
&\le \sup_{z \in [B, 2M^*)}\prob_{(-\beta/2,z)}(\tau_2(B) \wedge \tau_2(2M^*) > \log (2M^* B^{-1}))\\
&\hspace{7cm} + \sup_{z \in [B, 2M^*)} \prob_{(0,z)}\left(\tau_1(-\beta/2) \ge T(M^*)\right)
\end{split}
\end{equation}
where we have used the strong Markov property in the last step.
By \eqref{timelength},
$$
 \sup_{z \in [B, 2M^*)}\prob_{(-\beta/2,z)}\Big(\tau_2(B) \wedge \tau_2(2M^*) > \log (2M^* B^{-1})\Big) \le 1-\prob(\mathcal{E}(M^*)).
$$
By using $S(0) + \sqrt{2}W(t) - \beta t/2 \ge S(t)$ for $t \le \tau_1(-\beta/2)$ and $Q_1(t) \le S(t) \le Q_2(t)$ for $t \ge 0$,
\begin{multline*}
\sup_{z \in [B, 2M^*)} \prob_{(0,z)}\left(\tau_1(-\beta/2) \ge T(M^*)\right) \le \prob\left(\inf_{t \le T(M^*)}\left(\sqrt{2}W(t) - \beta t/2 \right) \ge  -(2M^* + \beta/2)\right)\\
\le \prob\left(\sqrt{2} W(T(M^*)) \ge \beta T(M^*)/2 - (2M^* + \beta/2)\right) \le \prob(\mathcal{E}(M^*))/2.
\end{multline*}
Using these bounds in \eqref{downadd1}, we obtain
\begin{equation}\label{timefin}
\sup_{z \in [B, 2M^*)} \prob_{(0,z)}\left(\tau_2(B) \wedge \tau_2(2M^*) > T(M^*) + \log (2M^* B^{-1})\right) \le 1 - \frac{\prob(\mathcal{E}(M^*))}{2} <1.
\end{equation}
Thus, using the strong Markov property and \eqref{timefin}, we obtain for any $k \ge 0$,
\begin{align}\label{down11}
\prob_{(0,y)}(\alpha^*_{2k+1} - \alpha^*_{2k} > n\left(T(M^*) + \log (2M^* B^{-1})\right)) \le \left(1 - \frac{\prob(\mathcal{E}(M^*))}{2}\right)^n.
\end{align}
Furthermore, by \eqref{down} we have constants $C_1$ and $C_2$, such that for $k \ge 1$ and for all $t \ge t^*$,
\begin{align}\label{down2}
\prob_{(0,y)}(\alpha^*_{2k} - \alpha^*_{2k-1}>t) \le C_1\exp(-C_2t^{1/5}).
\end{align}
Writing $\alpha_1 = \sum_{j=0}^{2\mathcal{N}'} (\alpha^*_{j+1} - \alpha^*_j)$ and using \eqref{down11} and \eqref{down2}, we get positive constants $C, C', C''$ and $t_{\alpha}^{(2)}>0$, depending on $\beta, B, M^*$, such that for all $t \ge t_{\alpha}^{(2)}$,
\begin{align*}
\prob_{(0,y)}(\alpha_1>t) &\le \prob(\mathcal{N}'>n) + \prob\Big(\sum_{j=0}^{2n} (\alpha^*_{j+1} - \alpha^*_j) > t\Big) \\
&\le \e^{-Cn} + C'n \e^{-C(t/n)^{1/5}} \le C'\e^{-C''t^{1/6}},
\end{align*}
where the last step is obtained by taking $n = \lfloor t^{1/6} \rfloor$.
\end{proof}

\subsection{Up-crossings of \texorpdfstring{$\mathbf{Q_2}$}{Q2}}\label{ssec:up}

In this subsection, we will prove tail-asymptotics for the distribution of $\alpha_2-\alpha_1$ as stated in Lemma~\ref{lem:alphauptail}.
The proof consists of the following two major parts: (i) First we establish in Lemma~\ref{lem:alpha} the tail probability of the hitting time of $Q_2$ to level $2B$ starting below level $B$ when $Q_1(0)$ is not too small.
(ii) Then in Lemma~\ref{qonealphaone} we show that at time $\alpha_1$, $Q_1(\alpha_1)$ cannot be too small.
Lemmas~\ref{lem:alpha} and~\ref{qonealphaone} are combined to prove Lemma~\ref{lem:alphauptail}.
\begin{lemma}\label{lem:alpha}
For any fixed $B>0$ and $M > 8B + 6\beta$, there exists $c_\alpha^{(2)}>0$ (depending on $M,B,\beta$) such that for all $t\geq 9$,
$$\sup_{x \in [-M/2,0], \ y \in (0,B]}\prob_{(x, y)}(\tau_2(2B)>t)\leq \exp(-c_\alpha^{(2)}\sqrt{t}).$$
\end{lemma}

To prove Lemma~\ref{lem:alpha} we split the time interval $[0, \tau_2(2B)]$ into subintervals
using stopping times of $Q_1$.
Depending on the local dynamics in each such subinterval we bound the diffusion process by more tractable diffusion processes, which are then used to obtain probability bounds.
The details of the proof of Lemma~\ref{lem:alpha} are given in Appendix~\ref{app:second}. 

As mentioned above, the next lemma gives a tail estimate on the distribution of $Q_1(\alpha_1)$.
\begin{lemma}\label{qonealphaone}
Fix $(Q_1(0), Q_2(0))=(x,y)$ with $x \le 0$, $y>0$. 
Recall the constant $t_{\alpha}^{(1)}$ obtained in Lemma~\ref{lem:alphadowntail}. There exist constants $C_1, C_2>0$ possibly depending on $(x,y)$, $B$, and $\beta$, such that for all $A \ge \max\{8 \beta t_{\alpha}^{(1)}, -4x\}$,
$$
\prob_{(x,y)}(Q_1(\alpha_1) < -A) \le C_1 \e^{-C_2 A^{1/6}}.
$$
\end{lemma}
\begin{proof}
In the proof, $C, C'$ will denote generic positive constants depending on $\beta, x, y$ whose values change from line to line. Observe that for $t>0$,
$$
Q_1(t) \ge Q_1(0) + \sqrt{2}W(t) - \beta t - L^*(t)
$$
where $L^*(t) = \sup_{s \le t}(Q_1(0) + \sqrt{2}W(s) - \beta s)^+$.
Thus, for any $A \ge \max\{8 \beta t_{\alpha}^{(1)}, -4x\}$,
\begin{equation}\label{Qal1}
\begin{split}
&\prob_{(x,y)}(Q_1(\alpha_1) < -A) \le \prob_{(x,y)}(\alpha_1 > A/(8\beta)) + \prob_{(x,y)}\left(\inf_{s \le A/(8\beta)}Q_1(s) < -A\right)\\
 &\le \prob_{(x,y)}(\alpha_1 > A/(8\beta)) + \prob_{(x,y)}\left(\inf_{s \le A/(8\beta)}\left(Q_1(0) + \sqrt{2}W(s) - \beta s - L^*(s)\right) < -A\right)\\
& \le \prob_{(x,y)}(\alpha_1 > A/(8\beta)) + \prob_{(x,y)}\left(L^*(A/(8\beta)) > A/2\right) \\
&\hspace{4cm}+ \prob_{(x,y)}\left(\inf_{s \le A/(8\beta)}\left(\sqrt{2}W(s) - \beta s\right) < -A/2 - x\right)\\
& \le \prob_{(x,y)}(\alpha_1 > A/(8\beta)) + \prob_{(x,y)}\left(\sup_{s \le A/(8\beta)}(\sqrt{2}W(s) - \beta s) > A/2\right)\\
  &\hspace{4cm}+ \prob_{(x,y)}\left(\inf_{s \le A/(8\beta)}\left(\sqrt{2}W(s) - \beta s\right) < -A/4\right).
\end{split}
\end{equation}
By Lemma \ref{lem:alphadowntail},
$$
\prob_{(x,y)}(\alpha_1 > A/(8\beta)) \le C\e^{-C'A^{1/6}}.
$$ 
Using the fact that the scale function (see \cite[V.46]{RW2000_2}) for $\sqrt{2}W(t) - \beta t$ is $s(z)= \exp(\beta z)$,
$$
\prob_{(x,y)}\left(\sup_{s \le A/(8\beta)}(\sqrt{2}W(s) - \beta s) > A/2\right) \le \prob_{(x,y)}\left(\sup_{s < \infty}(\sqrt{2}W(s) - \beta s) > A/2\right) = \e^{-\beta A/2}.
$$
Moreover, by standard estimates on normal distribution functions,
$$
\prob_{(x,y)}\left(\inf_{s \le A/(8\beta)}\left(\sqrt{2}W(s) - \beta s\right) < -A/4\right) \le \prob_{(x,y)}\left(\inf_{s \le A/(8\beta)}\left(\sqrt{2}W(s)\right) < -A/8\right) \le C\e^{-C'A}.
$$
Using the above bounds in \eqref{Qal1}, we obtain
$$
\prob_{(x,y)}(Q_1(\alpha_1) < -A) \le C\e^{-C'A^{1/6}}
$$
for any $A \ge \max\{8 \beta t_{\alpha}^{(1)}, -4x\}$, proving the lemma.
\end{proof}

\begin{proof}[Proof of Lemma~\ref{lem:alphauptail}]
In the proof, $C, C'$ will denote generic positive constants depending on $\beta, x, y$ whose values change from line to line. Fix $M > 8B + 6\beta + 2$. Take $t_{\alpha}' = 4\max\{9, 8 \beta t_{\alpha}^{(1)}, -4x,M\}$. Then for $t \ge t_{\alpha}'$,
\begin{multline}\label{upcross}
\prob_{(x,y)}(\alpha_2 - \alpha_1 >t) \le \prob_{(x,y)}(Q_1(\alpha_1) < -t/4) + \sup_{u \in[-t/4, -M/2], v>0}\prob_{(u,v)}(\tau_1(-M/2) > t/2)\\
 + \sup_{u \in [-M/2, 0], v \in (0,B]}\prob_{(u, v)}(\tau_2(2B) > t/2).
\end{multline}
By Lemma \ref{qonealphaone},
$$
\prob_{(x,y)}(Q_1(\alpha_1) < -t/4) \le C\e^{-C't^{1/6}}.
$$
Moreover,
\begin{multline*}
\sup_{u \in[-t/4, -M/2], v>0}\prob_{(u,v)}(\tau_1(-M/2) > t/2) \le \prob\left(-\frac{t}{4} + \sqrt{2}W(t/2) +\left(\frac{M}{2} - \beta\right)\frac{t}{2} < - \frac{M}{2}\right)\\
\le \prob\left(\sqrt{2}W(t/2) < - \frac{t}{4}\right) \le C\e^{-C't}.
\end{multline*}
By Lemma \ref{lem:alpha},
$$
\sup_{u \in [-M/2, 0], v \in (0,B]}\prob_{(u, v)}(\tau_2(2B) > t/2) \le \e^{-C'\sqrt{t}}.
$$
Using these bounds in \eqref{upcross}, we obtain for all $t \ge t_{\alpha}'$,
$$
\prob_{(x,y)}(\alpha_2 - \alpha_1 >t) \le C\e^{-C't^{1/6}}
$$
proving the lemma.
\end{proof}

\section{Analysis of fluctuations within a renewal cycle}\label{sec:fluc}
\begin{figure}
\begin{center}
\includegraphics[scale=.05]{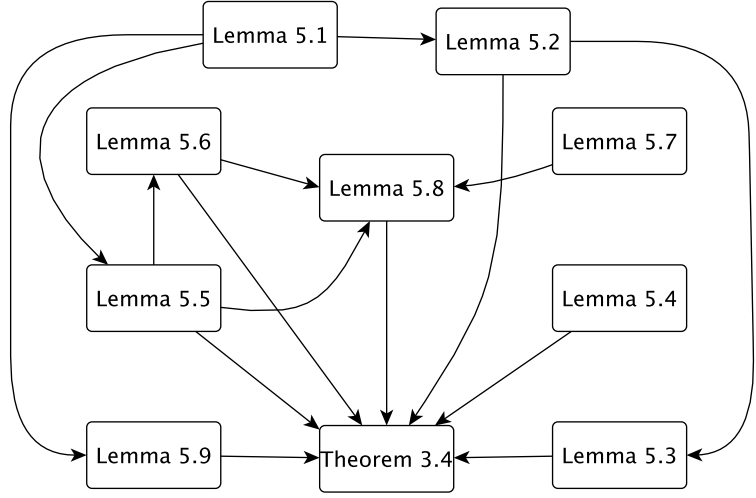}
\end{center}
\caption{Interdependence of various lemmas in the proof of Theorem~\ref{th:excren}. \label{fig:sec5}}
\end{figure}

In this section we prove Theorem~\ref{th:excren}.
Specifically,  we derive sharp estimates for the fluctuations of excursions of $Q_1$ and $Q_2$ between two successive regeneration times defined in \eqref{rendef}.
This will eventually furnish tail estimates for the stationary distribution of $Q_1$ and $Q_2$ and the scaling of extrema in large time intervals that are described in Theorems~\ref{th:statail} and~\ref{th:lil}.
First we state and prove Lemmas~\ref{lem:Q2hit1} -- \ref{lowqone}, which provide all the necessary results for proving Theorem~\ref{th:excren} at the end of this section.
For ease of understanding, in Figure~\ref{fig:sec5} we sketch the interdependence of various lemmas in this section.
\\

Denote the Brownian motion with drift $b$ and and its corresponding reflected analogue by
\begin{align*}
W^{(b)}(t) &:= \sqrt{2} W(t) + bt,\\
W^{(b)}_R(t) &:= \sqrt{2} W(t) + bt - \sup_{s \le t}\left(\sqrt{2} W(s) + bs\right),
\end{align*}
where $W$ denotes the standard Brownian motion.
Also, denote the local time of the reflected Brownian motion $W^{(b)}_R$ and its hitting time of level $z$ by $L^{(b)}$ and $\tau^{(b)}(z)$ respectively.

\begin{lemma}\label{lem:Q2hit1}
There exist positive constants $C_1, C_2>0$ that do not depend on $\beta$ such that
\begin{equation*}
\prob_{(0, y +\beta)}\left(\tau_1\left(-\frac{\beta}{2}\right) \le \tau_2\left(\frac{y}{2} + \beta\right)\right) \le C_1 \e^{-C_2\beta y}
\end{equation*}
for $y \ge \frac{1}{4\beta}$ if $\beta \ge 1$ and $y \ge \frac{64}{\beta}\log \frac{1}{\beta}$ if $\beta <1$.
\end{lemma}
\begin{proof}
From the evolution equation of $Q_1$ in~\eqref{eq:diffusionjsq}, note that for $y>0$, $W^{(y/2)}_R$ can be constructed on the same probability space as $(Q_1, Q_2)$, such that starting from $(Q_1(0), Q_2(0))= (0, y +\beta)$, almost surely $Q_1(t) \ge W^{(y/2)}_R(t)$ for all $t \le  \tau_2\left(\frac{y}{2} + \beta\right)$. 
The scale function $s$ for $W^{(y/2)}_R(t)$ is obtained by solving the equation $\frac{y}{2}s'(z) + s"(z)=0$ (see \cite[V.46]{RW2000_2}) and one candidate is
\begin{equation}\label{eq:scalf}
s(z) = \frac{2}{y} \left(1- \e^{-yz/2}\right).
\end{equation}
We will estimate the time taken by $W^{(y/2)}_R$ to hit the level $-\beta/2$. Define stopping times for the process $W^{(y/2)}$ as follows: For $i \ge 0$
\begin{align*}
\gamma_{i+1} &= \inf\big\{t \ge \gamma_{i} : W^{(y/2)}(t) - W^{(y/2)}(\gamma_i) \text{ hits } \beta/4 \text{ or } -\beta/4\big\},
\end{align*}
with the convention that $\gamma_0=0$.
From the explicit form of the scale function $s$ in~\eqref{eq:scalf}, observe that for $i \ge 0$,
\begin{equation}\label{scalebound}
\prob\left(W^{(y/2)}(\gamma_{i+1}) - W^{(y/2)}(\gamma_i) = -\beta/4\right) = \frac{1-\e^{-\beta y/8}}{\e^{\beta y/8} - \e^{-\beta y/8}} \le \e^{-\beta y /8}.
\end{equation}
Define
$$
\mathcal{N} := \inf\big\{ i \ge 1: W^{(y/2)}(\gamma_{i+1}) - W^{(y/2)}(\gamma_i) = -\beta/4\big\}.
$$
Then for any $n \ge 1$, by \eqref{scalebound}, $\prob(\mathcal{N} \le n) \le n \e^{-\beta y /8}$. Note that for $t < \gamma_{\mathcal{N}}$, $W^{(y/2)}_R(t) > - \beta/2$. Thus, $\tau^{(y/2)}\left(-\frac{\beta}{2}\right) \ge \gamma_{\mathcal{N}}$. 
Consequently,
$$
L^{(y/2)}\left(\tau^{(y/2)}\left(-\frac{\beta}{2}\right)\right) \ge \sup_{t \le \gamma_{\mathcal{N}}}\left(W^{(y/2)}(t)\right) \ge \mathcal{N}\beta/4.
$$
Therefore, for any $n \ge 1$,
\begin{equation}\label{locbound}
\prob\left(L^{(y/2)}\left(\tau^{(y/2)}\left(-\frac{\beta}{2}\right)\right) \le n \beta \right) \le \prob\left(\mathcal{N} \le 4n\right) \le 4n \e^{-\beta y /8}.
\end{equation}
Further, on the event $\left[ \tau_1\left(-\frac{\beta}{2}\right) \le \tau_2\left(\frac{y}{2} + \beta\right)\right]$, $\tau_1\left(-\frac{\beta}{2}\right) \ge \tau^{(y/2)}\left(-\frac{\beta}{2}\right)$. Therefore, for $n \ge 1$,
\begin{equation}\label{tau1}
\begin{split}
&\prob_{(0, y+ \beta)}\left(\tau_1\left(-\frac{\beta}{2}\right) \le n\beta/y, \tau_1\left(-\frac{\beta}{2}\right) \le \tau_2\left(\frac{y}{2} + \beta\right)\right)
\le \prob\left(\tau^{(y/2)}\left(-\frac{\beta}{2}\right) \le n\beta/y\right)\\
&\le \prob\left(\tau^{(y/2)}\left(-\frac{\beta}{2}\right) \le n\beta/y, L^{(y/2)}\left(\tau^{(y/2)}\left(-\frac{\beta}{2}\right)\right) > n \beta\right)
 + \prob\left(L^{(y/2)}\left(\tau^{(y/2)}\left(-\frac{\beta}{2}\right)\right) \le n \beta\right).
\end{split}
\end{equation}
An upper bound for the second probability in the right side of~\eqref{tau1} has been obtained in \eqref{locbound}. To estimate the first probability, observe that
\begin{equation}\label{tau2}
\begin{split}
&\prob\left(\tau^{(y/2)}\left(-\frac{\beta}{2}\right) \le n\beta/y, L^{(y/2)}\left(\tau^{(y/2)}\left(-\frac{\beta}{2}\right)\right) > n \beta\right)\\
&\hspace{4cm}\le \prob\left(\sup_{t \le n\beta/y}\left(\sqrt{2}W(t) + yt/2\right) > n\beta\right)\\
& \hspace{4cm}\le \prob\left(\sup_{t \le n\beta/y}\sqrt{2}W(t) > n\beta/2\right) \le \frac{4}{\sqrt{\pi n\beta y}} \e^{-n\beta y/16}.
\end{split}
\end{equation}
Using \eqref{locbound} and \eqref{tau2} in \eqref{tau1}, we obtain
\begin{multline}\label{imp1}
\prob_{(0, y+ \beta)}\left(\tau_1\left(-\frac{\beta}{2}\right) \le n\beta/y, \tau_1\left(-\frac{\beta}{2}\right) \le \tau_2\left(\frac{y}{2} + \beta\right)\right) \le 4n \e^{-\beta y /8} + \frac{4}{\sqrt{\pi n\beta y}} \e^{-n\beta y/16},
\end{multline}
where an appropriate choice of $n \ge 1$ (depending on $y$ and $\beta$) will be made later. Now, we want to estimate the probability $\prob_{(0, y+ \beta)}\left(n\beta/y < \tau_1\left(-\frac{\beta}{2}\right) \le \tau_2\left(\frac{y}{2} + \beta\right)\right)$. Towards this end, recall that $S(t) = Q_1(t) + Q_2(t)$ has the representation
$$
S(t) = S(0) + \sqrt{2}W(t) - \beta t + \int_0^t (-Q_1(s))ds.
$$
Thus, for $t \le \tau_1\left(-\frac{\beta}{2}\right)$,
$$
S(t) \le S(0) + \sqrt{2}W(t) - \frac{\beta}{2} t.
$$
Therefore, if $n$ is chosen such that $y \le \sqrt{n} \beta/4$,
\begin{equation}\label{imp2}
\begin{split}
&\prob_{(0, y+ \beta)}\left(n\beta/y < \tau_1\left(-\frac{\beta}{2}\right) \le \tau_2\left(\frac{y}{2} + \beta\right)\right)\\
 &\le \prob\left(y+\beta + \sqrt{2}W(t) - \frac{\beta}{2} t \ge y/2 + \beta/2, \text{ for all } t \le n\beta/y\right)\\
 &\le \prob\left(\sqrt{2}W(n\beta/y) - \frac{n\beta^2}{2y} \ge -y/2 - \beta/2\right)\\
& \le \prob\left(\sqrt{2}W(n\beta/y) \ge \frac{n\beta^2}{8y}\right) \qquad \text{since } y \le \sqrt{n} \beta/4\\
 &\le \frac{8\sqrt{y}}{\sqrt{\pi n} \beta^{3/2}} \e^{-\frac{n\beta^3}{256 y}}.
\end{split}
\end{equation}
From \eqref{imp1} and \eqref{imp2} we obtain
\begin{multline}\label{imp3}
\prob_{(0, y+ \beta)}\left(\tau_1\left(-\frac{\beta}{2}\right) \le \tau_2\left(\frac{y}{2} + \beta\right)\right) \le 4n \e^{-\beta y /8} + \frac{4}{\sqrt{\pi n\beta y}} \e^{-n\beta y/16} + \frac{8\sqrt{y}}{\sqrt{\pi n} \beta^{3/2}} \e^{-\frac{n\beta^3}{256 y}}.
\end{multline}
Now, if $\beta \ge 1$, choose $n = 16 y^2\beta^2$. Then, clearly $y \le \sqrt{n} \beta/4$. With this choice of $n$, the above expression yields the following bound:
\begin{multline}\label{imp4}
\prob_{(0, y+ \beta)}\left(\tau_1\left(-\frac{\beta}{2}\right) \le \tau_2\left(\frac{y}{2} + \beta\right)\right) \le 64(\beta y)^2 \e^{-\beta y /8} + \frac{1}{\sqrt{\pi}(\beta y)^{3/2}} \e^{-(\beta y)^3} + \frac{2}{\sqrt{\pi \beta y}} \e^{-\frac{\beta y}{16}}
\end{multline}
for $y \ge \frac{1}{4\beta}$ (this ensures $n \ge 1$).\\

If $\beta < 1$, choose $n=y^4$. 
Then $y \le \sqrt{n} \beta/4$ is satisfied if $y \ge 4/\beta$. Some routine calculations reveal that for $y \ge 4/\beta$ the second and third terms appearing on the right side of~\eqref{imp3} can be estimated by,
$$
\frac{4}{\sqrt{\pi n\beta y}} \e^{-n\beta y/16} \le \frac{1}{8\sqrt{\pi}} \e^{-(\beta y)^5/16}
$$
and
$$
\frac{8\sqrt{y}}{\sqrt{\pi n} \beta^{3/2}} \e^{-\frac{n\beta^3}{256 y}} \le \frac{1}{\sqrt{\pi}}\e^{-(\beta y)^3/256}.
$$
To estimate the first term on the right side of~\eqref{imp3}, rewrite it as
$$
4n \e^{-\beta y /8} = \left[4(\beta y)^4 \e^{-(\beta y)/16}\right] \left[\beta^{-4}\e^{-(\beta y)/16}\right].
$$
Observe that $\beta^{-4}\e^{-(\beta y)/16} \le 1$ for $y \ge \frac{64}{\beta}\log \frac{1}{\beta}$. Therefore, for $\beta<1$ and $y \ge \frac{64}{\beta}\log \frac{1}{\beta}$, we have the following bound:
\begin{equation}\label{imp5}
\prob_{(0, y+ \beta)}\left(\tau_1\left(-\frac{\beta}{2}\right) \le \tau_2\left(\frac{y}{2} + \beta\right)\right) \le \frac{1}{8\sqrt{\pi}} \e^{-(\beta y)^5/16} + \frac{1}{\sqrt{\pi}}\e^{-(\beta y)^3/256} + 4(\beta y)^4 \e^{-(\beta y)/16}.
\end{equation}
The lemma follows from \eqref{imp4} and \eqref{imp5}.
\end{proof}
The above lemma can be used to deduce the following hitting time estimate for $Q_2$.
\begin{lemma}\label{Q2gebeta1}
There exist constants $\widetilde{C}_1, \widetilde{C}_2 >0$ that do not depend on $\beta$ such that
\begin{equation*}
\prob_{(0, y+ \beta)}\left(\tau_2\left(2y + \beta\right) \le \tau_2\left(\frac{y}{2} + \beta\right)\right) \le \widetilde{C}_1 \e^{-\widetilde{C}_2 \beta y}
\end{equation*}
for $y \ge \frac{1}{4\beta}$ if $\beta \ge 1$ and $y \ge \frac{64}{\beta}\log \frac{1}{\beta}$ if $\beta <1$.
\end{lemma}
\begin{proof}
We can write for any $y>0$,
\begin{multline}\label{Q21}
\prob_{(0, y+ \beta)}\left(\tau_2\left(2y + \beta\right) \le \tau_2\left(\frac{y}{2} + \beta\right)\right)\\
\le \prob_{(0, y+ \beta)}\left(\tau_1\left(-\frac{\beta}{2}\right) \le \tau_2\left(\frac{y}{2} + \beta\right)\right) + \prob_{(0, y+ \beta)}\left(\tau_2\left(2y + \beta\right) < \tau_1\left(-\frac{\beta}{2}\right)\right).
\end{multline}
By Lemma \ref{lem:Q2hit1},
\begin{equation}\label{Q22}
\prob_{(0, y+ \beta)}\left(\tau_1\left(-\frac{\beta}{2}\right) \le \tau_2\left(\frac{y}{2} + \beta\right)\right) \le C_1\e^{-C_2\beta y}
\end{equation}
for $y \ge \frac{1}{4\beta}$ if $\beta \ge 1$ and $y \ge \frac{64}{\beta}\log \frac{1}{\beta}$ if $\beta <1$.
To estimate the second probability in \eqref{Q21}, recall that for $t \le \tau_1\left(-\frac{\beta}{2}\right)$, $S(t) = Q_1(t) + Q_2(t)$ satisfies
$$
S(t) \le S(0) + \sqrt{2}W(t) - \frac{\beta}{2} t.
$$
Therefore,
\begin{align}\label{Q23}
\prob_{(0, y+ \beta)}\left(\tau_2\left(2y + \beta\right) < \tau_1\left(-\frac{\beta}{2}\right)\right)
\le \prob\left( \sup_{t< \infty}\left(\sqrt{2}W(t) - \frac{\beta}{2} t\right) \ge y\right) = \e^{-\frac{\beta y}{2}}
\end{align}
for $y >0$. The first inequality above follows from the fact that points of time where $Q_2$ increases are precisely those where $Q_1$ equals zero: hence $Q_1(\tau_2(2y+\beta))=0$.

The lemma now follows by using \eqref{Q22} and \eqref{Q23} in \eqref{Q21}.
\end{proof}
The above estimate can be strengthened to the following tail estimate which will be used to study fluctuations of $Q_2$ between successive regeneration times.
\begin{lemma}\label{Q2gebeta2}
Recall the constants $\widetilde{C}_1, \widetilde{C}_2$ in the statement of Lemma \ref{Q2gebeta1}. There exist constants $C^*_1, C^*_2 >0$ that do not depend on $\beta$ such that
\begin{equation*}
\prob_{(0, y+ \beta)}\left(\tau_2\left(2y + \beta\right) \le \tau_2\left(y_0 + \beta\right)\right) \le C^*_1 \e^{-C^*_2 \beta y}
\end{equation*}
for all $y \ge y_0$, where $y_0 = \max\left\lbrace \frac{1}{4\beta}, \frac{\log(4\widetilde{C}_1)}{\widetilde{C}_2\beta}\right\rbrace$ if $\beta \ge 1$ and $y_0 = \max\left\lbrace\frac{64}{\beta}\log \frac{1}{\beta}, \frac{\log(4\widetilde{C}_1)}{\widetilde{C}_2\beta}\right\rbrace$ if $\beta <1$.
\end{lemma}
\begin{proof}
Define stopping times:
\begin{align*}
T_{2k+1}&=\inf\left\lbrace t \ge T_{2k}: Q_2(t)=2y + \beta \text{ or } Q_2(t)=\frac{y}{2} + \beta \text{ or } Q_2(t)=y_0 + \beta \right\rbrace;\\
T_{2k+2}&=\inf\left\lbrace t \ge T_{2k+1}: Q_2(t)=y + \beta \text{ or } Q_2(t)=y_0 + \beta \right\rbrace,
\end{align*}
for $k \ge 0$, with the convention that $T_0=0$. Let
$$
\mathcal{N}^{0} = \inf\{ k \ge 1: Q_2\left(T_{2k}\right) = y_0 + \beta\}.
$$
Define $\hat{Q}_2(t) = \log_2(Q_2(t) - \beta)$. By Lemma \ref{Q2gebeta1} and our choice of $y_0$, for any $z \ge \log_2(y_0)$,
\begin{multline*}
\prob(\hat{Q}_2 \text{ hits } z+1 \text{ before } z-1 \ \vert \ \hat{Q}_2(0) = z, Q_1(0)=0)\\
=\prob_{(0,2^z + \beta)}(Q_2 \text{ hits } 2^{z+1}+\beta \text{ before } 2^{z-1} + \beta) \le 1/4.
\end{multline*}
Thus, $\hat{Q}_2$ starting from any $z \ge \log_2(y_0)$ and observed at the stopping times where the increments are $\pm1$ until the first time it crosses the level $\log_2(y_0)$ (i.e., strictly less than $\log_2(y_0)$) is stochastically dominated by a random walk $\left(\mathcal{S}_n\right)_{n \ge 0}$ where 
$$\prob(S_{n+1} - S_n =1) = 1- \prob(S_{n+1} - S_n =-1) = 1/4.$$ Therefore,
\begin{multline*}
\sup_{z \ge \log_2(y_0)}\prob(\hat{Q}_2 \text{ hits } z+1 \text{ before it crosses } \log_2(y_0) \ \vert \ \hat{Q}_2(0) = z, Q_1(0)=0)\\
\le \sup_{z \ge \log_2(y_0)}\prob(S_n \text{ hits } z+1 \ \vert \ S_0=z) = p^{(S)} <1.
\end{multline*}
which, in turn, implies that for any $y \ge y_0$,
\begin{align*}
\prob_{(0, y+ \beta)}\left(\tau_2\left(2y + \beta\right) \le \tau_2\left(y_0 + \beta\right)\right) \le p^{(S)} <1.
\end{align*}
Thus, for any $k \ge 1$,
\begin{equation}\label{Nzero}
\prob_{(0, y+ \beta)}\left(\mathcal{N}^{0} \ge k+1\right) \le (p^{(S)})^k.
\end{equation}
Finally, for any $y \ge y_0$,
\begin{align*}
&\prob_{(0, y+ \beta)}\left(\tau_2\left(2y + \beta\right) \le \tau_2\left(y_0 + \beta\right)\right) \\
&=\prob_{(0, y+ \beta)}\left(\sup_{0 \le t \le T_{2\mathcal{N}^0}}Q_2(t) > 2y + \beta\right)
\le \sum_{k=1}^{\infty}\prob_{(0, y+ \beta)}\left(\sup_{T_{2k-2} \le t \le T_{2k}}Q_2(t) > 2y + \beta, \mathcal{N}^0 \ge k\right)\\
 &=  \sum_{k=1}^{\infty}\mathbb{E}_{(0, y+ \beta)}\mathbb{I}(\mathcal{N}^0 \ge k) \prob_{(0, y+ \beta)}\left(\tau_2\left(2y + \beta\right) \le \tau_2\left(\frac{y}{2} + \beta\right)\right), \ \text{by strong Markov property at } T_{2k-2}\\
& \le \prob_{(0, y+ \beta)}\left(\tau_2\left(2y + \beta\right) \le \tau_2\left(\frac{y}{2} + \beta\right)\right) \sum_{k=1}^{\infty} (p^{(S)})^{k-1}, \qquad \text{by } \eqref{Nzero}\\
 &\le (1-p^{(S)})^{-1} \widetilde{C}_1 \e^{-\widetilde{C}_2 \beta y},  \qquad \text{by Lemma } \ref{Q2gebeta1},
\end{align*}
which completes the proof of the lemma.
\end{proof}
The lower bound on the tail probabilities is achieved for all $\beta>0$ in the following lemma.
\begin{lemma}\label{Q2lb}
For any $\beta>0$ and any $B>0$,
\begin{equation*}
\prob_{(0,2B)}\left(\tau_2(y) < \tau_2(B)\right) \ge (1-\e^{-\beta B})\e^{-\beta(y-2B)}
\end{equation*}
for all $y \ge 2B$.
\end{lemma}
\begin{proof}
Note that $Q_2(t) \ge Q_1(t) + Q_2(t) = S(t)$ for all $t\ge 0$. Further, recall that
$$
S(t) = S(0) + \sqrt{2}W(t) - \beta t + \int_0^t(-Q_1(s))ds \ge S(0) + \sqrt{2}W(t) - \beta t, \ t \ge 0.
$$
Therefore, for all $y \ge 2B$,
\begin{align*}
\prob_{(0,2B)}\left(\tau_2(y) < \tau_2(B)\right) &\ge \prob_{(0,2B)}\left(S(t) \text{ hits level } y \text{ before level } B\right)\\
&\ge \prob\left(2B + \sqrt{2}W(t) - \beta t \text{ hits level } y \text{ before level }B\right)\\
&=\prob\left(\sqrt{2}W(t) - \beta t \text{ hits level } y-2B \text{ before level } -B\right)\\
&= \frac{1-\e^{-\beta B}}{\e^{\beta(y-2B)} - \e^{-\beta B}}, \qquad \text{by scale function arguments}\\
&\ge (1-\e^{-\beta B})\e^{-\beta(y-2B)},
\end{align*}
proving the lemma.
\end{proof}
Now, we will study fluctuations of $Q_1$ within one renewal cycle. 
Recall 
$l_0(\beta)$ from~\eqref{lzerodef} and the notation
$$
\sigma(t) =\inf\{ s \ge t: Q_1(s)=0\}, \ \ t \ge 0.
$$
\begin{lemma}\label{unifprob}
There exist constants $R_1>0$ not depending on $\beta$ and $p^{**}(\beta) \in (0,1)$ such that for all $R \ge R_1$,
\begin{equation}
\sup_{y \ge Rl_0(\beta)} \prob_{(0,y)}\left(\tau_1(-\beta) < \tau_2(Rl_0(\beta))\right) = p^*(\beta,R)\leq p^{**}(\beta).
\end{equation}
\end{lemma}
\begin{proof}
In the proof $C,C',C_1, C_2, \dots$ will denote generic positive constants not depending on $\beta, R$ whose values might change from line to line. For any $y \ge Rl_0(\beta) - \beta$,
\begin{multline}\label{R11}
\prob_{(0,y+\beta)}\left(\tau_1(-\beta) < \sigma\left(\tau_2\left(\frac{y}{2} + \beta\right)\right)\right)
 \le \prob_{(0, y+\beta)} \left(\tau_1(-\beta/2) \le \tau_2\left(\frac{y}{2} + \beta\right)\right)\\
  + \prob_{(0,y+\beta)}\left(\tau_2\left(\frac{y}{2} + \beta\right) < \tau_1(-\beta/2) < \tau_1(-\beta) < \sigma\left(\tau_2\left(\frac{y}{2} + \beta\right)\right)\right) \\
  \le \prob_{(0, y+\beta)} \left(\tau_1(-\beta/2) \le \tau_2\left(\frac{y}{2} + \beta\right)\right) + \sup_{x \in [-\beta/2, 0]}\prob_{\left(x,\frac{y}{2} + \beta\right)}\left(\tau_1(-\beta) < \tau_1(0)\right)
\end{multline}
where the last step is a consequence of the strong Markov property applied at $\tau_2\left(\frac{y}{2} + \beta\right)$.
From Lemma \ref{lem:Q2hit1}, for $R \ge 65$,
\begin{equation}\label{Rone}
\prob_{(0, y+\beta)} \left(\tau_1(-\beta/2) \le \tau_2\left(\frac{y}{2} + \beta\right)\right) \le C_1\e^{-C_2\beta y}, \ \ y \ge Rl_0(\beta) - \beta.
\end{equation}
Now let us take the starting configuration to be $(Q_1(0), Q_2(0)) = \left(x,\frac{y}{2} + \beta\right)$ with $y \ge Rl_0(\beta) - \beta$ and $R \ge 5$.
In that case, since $(\dif/\dif t)Q_2(t)\geq -Q_2(t)$, therefore $Q_2(t)\geq (y/2 +\beta)/2$ for all $t \le \log 2$.
Consequently, for any $t \le \log 2$,
\begin{align*}
Q_1(t) &= Q_1(0) + \sqrt{2} W(t) - \beta t +
\int_0^t (- Q_1(s) + Q_2(s)) \dif s - L(t)\\
&\geq x + \sqrt{2} W(t) - \beta t +
\int_0^t Q_2(s) \dif s\\
&\ge x + \sqrt{2}W(t) + (y-2\beta)t/4 \ge x + \sqrt{2}W(t) + yt/8.
\end{align*}
Therefore,
\begin{equation}\label{R12}
\begin{split}
&\sup_{x \in [-\beta/2, 0]}\prob_{\left(x,\frac{y}{2} + \beta\right)}\left(\tau_1(-\beta) < \tau_1(0) \le \log 2\right)\\
&\hspace{2cm} \le \sup_{x \in [-\beta/2, 0]}\prob_{\left(x,\frac{y}{2} + \beta\right)}\left(x + \sqrt{2}W(t) + yt/8 \text{ hits } -\beta \text{ before } 0\right)\\
&\hspace{2cm}  \le \prob\left(\sqrt{2}W(t) + yt/8 \text{ hits } -\beta/2 \text{ before } \beta/2\right) \le \e^{-\beta y/16}
\end{split}
\end{equation}
where the last step follows from standard scale function arguments. Moreover, for $y \ge Rl_0(\beta) - \beta$ with $R \ge 65$,
\begin{equation}\label{R13}
\begin{split}
&\sup_{x \in [-\beta/2, 0]}\prob_{\left(x,\frac{y}{2} + \beta\right)}\left(\tau_1(0) > \log 2\right)\\
 &\le  \sup_{x \in [-\beta/2, 0]}\prob_{\left(x,\frac{y}{2} + \beta\right)}\left(\sup_{t \le \log 2}(x + \sqrt{2}W(t) + yt/8) < 0\right)\\
&\le\prob\left(\sup_{t \le \log 2}( \sqrt{2}W(t) + yt/8) < \beta/2\right)
 \le \prob\left(\sqrt{2}W(\log 2)  < -y/32\right)\\
 &\le \e^{-y^2/(4(32^2) \log 2)} \le \e^{-\beta y/(64 \log 2)}.
\end{split}
\end{equation}
Using \eqref{Rone}, \eqref{R12} and \eqref{R13} in \eqref{R11}, we obtain for $R \ge 65$, there exist positive constants $C,C'$ not depending on $\beta$ and $R$ such that for all $y \ge Rl_0(\beta) - \beta$
\begin{equation}\label{Rsigma}
\prob_{(0,y+\beta)}\left(\tau_1(-\beta) < \sigma\left(\tau_2\left(\frac{y}{2} + \beta\right)\right)\right) \le C\e^{-C'\beta y}.
\end{equation}
Now, for any $y \ge Rl_0(\beta) - \beta$, observe that the event  $[\tau_1(-\beta/2) \le \tau_2(Rl_0(\beta))]$ can be written as
\[\big[\tau_1(-\beta) \le \tau_2(Rl_0(\beta))\big] \subseteq \bigcup_{k=1}^{\left\lfloor\log_2\left(\frac{y}{Rl_0(\beta)-\beta}\right) + 2\right\rfloor}\Big[\sigma\left(\tau_2\left(\frac{y}{2^{k-1}} + \beta\right)\right) < \tau_1\left(-\beta\right) < \sigma\left(\tau_2\left(\frac{y}{2^k} + \beta\right)\right)\Big],\]
and therefore,
\begin{multline}\label{Rtwo}
\prob_{(0,y+\beta)}\left(\tau_1(-\beta) \le \tau_2(Rl_0(\beta))\right)\\
\le \sum_{k=1}^{\left\lfloor\log_2\left(\frac{y}{Rl_0(\beta)-\beta}\right) + 2\right\rfloor}\prob_{(0,y+\beta)}\left(\sigma\left(\tau_2\left(\frac{y}{2^{k-1}} + \beta\right)\right) < \tau_1(-\beta) < \sigma\left(\tau_2\left(\frac{y}{2^k} + \beta\right)\right)\right).
\end{multline}
Take any $R \ge 260$. By the strong Markov property, for each $k \le \left\lfloor\log_2\left(\frac{y}{Rl_0(\beta)-\beta}\right) + 2\right\rfloor$,
\begin{align*}
&\prob_{(0,y+\beta)}\left(\sigma\left(\tau_2\left(\frac{y}{2^{k-1}} + \beta\right)\right) < \tau_1(-\beta) < \sigma\left(\tau_2\left(\frac{y}{2^k} + \beta\right)\right)\right)\\
&\hspace{3cm} \le \sup_{z \in [y/2^k, y/2^{k-1}]} \prob_{(0,z+\beta)}\left(\tau_1(-\beta) < \sigma\left(\tau_2\left(\frac{y}{2^k} + \beta\right)\right)\right)\\
&\hspace{3cm}\le \sup_{z \in [y/2^k, y/2^{k-1}]} \prob_{(0,z+\beta)}\left(\tau_1(-\beta) < \sigma\left(\tau_2\left(\frac{z}{2} + \beta\right)\right)\right) \le C\e^{-C'\beta y/2^k},
\end{align*}
where the last inequality follows from \eqref{Rsigma} as for $k \le \left\lfloor\log_2\left(\frac{y}{Rl_0(\beta)-\beta}\right) + 2\right\rfloor$, $\frac{y}{2^k} \ge \frac{Rl_0(\beta) - \beta}{4} \ge \frac{R}{4}l_0(\beta) - \beta$ and $\frac{R}{4} \ge 65$.

Writing $p(\beta,R) = C_1\e^{-C_2\beta (Rl_0(\beta)-\beta)/4}$ and using the above bound in \eqref{Rtwo}, we obtain $R_1 > 0$ such that for any $R \ge R_1$ and any $y \ge Rl_0(\beta) - \beta$,
\begin{align}
\prob_{(0,y+\beta)}\left(\tau_1(-\beta) \le \tau_2(Rl_0(\beta))\right) &\le \sum_{k=1}^{\left\lfloor\log_2\left(\frac{y}{Rl_0(\beta)-\beta}\right) + 2\right\rfloor}C_1\e^{-C_2\beta y/2^k}\label{eq:lem5.5sum1}\\
& \le \sum_{k=0}^{\infty}p(\beta,R)^{2^k} \le \sum_{k=0}^{\infty}p(\beta,R_1)^{2^k} =: p^{**}(\beta) <1,\label{eq:lem5.5sum2}
\end{align}
where the second inequality can be seen as follows: For any $y\ge Rl_0(\beta) - \beta$, the last term in the sum in~\eqref{eq:lem5.5sum1} is bounded above by $p(\beta, R)$.
Also, starting from the last term and counting backwards in $k$, observe that each next term is the square of the previous term, which provides the $2^k$ in the exponent of $p(\beta, R)$ in~\eqref{eq:lem5.5sum2}.
Now, it is straightforward to see that for a fixed $\beta$ the first sum in \eqref{eq:lem5.5sum2} is a decreasing function in $R$, and is bounded away from 1 for all large enough $R$.
This proves the lemma.
\end{proof}
\begin{lemma}\label{alphaone}
There exists a constant $R_2>0$ not depending on $\beta$ such that for any $R \ge R_2$, there is a constant $C_2(\beta,R)>0$ (depending on $\beta, R$) satisfying
\begin{equation}
\sup_{z \in[-\beta, 0],y \ge 2R l_0(\beta)}\prob_{(z,y)}\left(\tau_1(-x) < \tau_2(Rl_0(\beta))\right) \le C_2(\beta,R)\e^{-(x-\beta)^2/2}, \ \ \text{ for all } x \ge \beta +1.
\end{equation}
\end{lemma}
\begin{proof}
Take any $R>0$. Let $(Q_1(0), Q_2(0))=(z,y)$ where $z \in [-\beta,0]$ and $y \ge 2Rl_0(\beta)$. Define the stopping times: $\sigma^{(0)} =0$ and for $k \ge 0$,
\begin{align*}
\sigma^{(2k+1)} &= \inf\{t \ge \sigma^{(2k)}: Q_1(t) = -\beta-1 \text{ or } Q_2(t) \le Rl_0(\beta)\},\\
\sigma^{(2k+2)} &= \inf\{t \ge \sigma^{(2k+1)}: Q_1(t) = -\beta \text{ or } Q_2(t) \le Rl_0(\beta)\}.
\end{align*}
Define $\mathcal{N}^{\sigma} = \inf\{n \ge 1: Q_2(\sigma^{(n)}) \le Rl_0(\beta)\}$. Observe that for any $z \in[-\beta, 0]$, by the strong Markov property, we obtain
\begin{equation}\label{eq:lem5.6-1}
\begin{split}
&\sup_{y \ge R l_0(\beta)}\prob_{(z,y)}\left(\tau_1(-\beta-1) < \tau_2(Rl_0(\beta))\right)\\
& \hspace{2cm}\le \sup_{y \ge R l_0(\beta)}\prob_{(z,y)}\left(\tau_1(0) < \tau_1(-\beta-1) < \tau_2(Rl_0(\beta))\right)\\
& \hspace{5cm} + \sup_{y \ge R l_0(\beta)}\prob_{(z,y)}\left(\tau_1(-\beta-1) < \tau_1(0) \wedge  \tau_2(Rl_0(\beta))\right)\\
&\hspace{2cm}\le \sup_{y \ge R l_0(\beta)}\prob_{(0,y)}\left(\tau_1(-\beta) < \tau_2(Rl_0(\beta))\right) \\
& \hspace{5cm}+ \sup_{y \ge R l_0(\beta)}\prob_{(z,y)}\left(\tau_1(-\beta-1) < \tau_1(0)\wedge  \tau_2(Rl_0(\beta))\right).
\end{split}
\end{equation}
By Lemma \ref{unifprob}, for large enough $R$,
\begin{equation}\label{lal}
\sup_{y \ge R l_0(\beta)}\prob_{(0,y)}\left(\tau_1(-\beta) < \tau_2(Rl_0(\beta))\right) \le p^{**}(\beta) <1.
\end{equation}
Further, observe that for $t \le \tau_1(0)\wedge  \tau_2(Rl_0(\beta))$, 
$$Q_1(t) \ge z + \sqrt{2}W(t) + (Rl_0(\beta) - \beta)t \ge -\beta + \sqrt{2}W(t) + (Rl_0(\beta) - \beta)t.$$ 
Therefore,
\begin{multline}\label{nil}
 \sup_{y \ge R l_0(\beta)}\prob_{(z,y)}\left(\tau_1(-\beta-1) < \tau_1(0)\wedge  \tau_2(Rl_0(\beta))\right)\\
 \le \prob(-\beta + \sqrt{2}W(t) + (Rl_0(\beta) - \beta)t \text{ hits } -\beta -1 \text{ before } 0) \le \e^{-(Rl_0(\beta)-\beta)}.
\end{multline}
Using~\eqref{lal} and~\eqref{nil} in~\eqref{eq:lem5.6-1}, we conclude that there is $R_2>0$ such that for all $R \ge R_2$,
\begin{equation}\label{holud}
\sup_{z \in[-\beta, 0],y \ge R l_0(\beta)}\prob_{(z,y)}\left(\tau_1(-\beta-1) < \tau_2(Rl_0(\beta))\right) \le p'(\beta,R)<1.
\end{equation}
Using \eqref{holud} and the strong Markov property, there exists a constant $C(\beta,R)>0$ depending on $\beta,R$ such that
\begin{equation}\label{Nexpo}
\sup_{z \in[-\beta, 0],y \ge 2R l_0(\beta)}\mathbb{E}_{(z,y)}(\mathcal{N}^{\sigma}) \  \le 2\sum_{n=0}^{\infty}\prob(\mathcal{N}^{\sigma} > 2n)  \ \le 2\sum_{n=0}^{\infty}p'(\beta, R)^n \ \le C(\beta,R) \ <\infty.
\end{equation}
For $(Q_1(0), Q_2(0)) = (-\beta-u, y)$ for any $u\ge 1,y > 0$, by \cite[Proposition 2.18]{Karatzas}, a process $Z$ can be constructed on the same probability space as $(Q_1,Q_2)$, such that $Q_1(t) + \beta \ge Z(t)$ for $t \le \tau_1(0)$, where $Z$ is an Ornstein-Uhlenbeck process which solves the SDE:
$$
dZ(t) = \sqrt{2}dW(t) - Z(t)dt, \ \ Z(0)=-u.
$$
The scale function for $Z$ is given by $s_Z(z) = \int_0^z\e^{w^2/2}dw$. From this observation and elementary estimates on $s_Z$, we have for any $x \ge \beta+u$,
\begin{multline}\label{OUcom}
\sup_{y > 0}\prob_{(-\beta-u, y)}\left(\tau_1(-x) < \tau_1(-\beta)\right) \le \prob\left(Z(t) \text{ hits } -x+\beta \text{ before } 0\right)\\
= \frac{s_Z(0) - s_Z(-u)}{s_Z(0)-s_Z(-x+\beta)}\le \sqrt{9\pi/2}\e^{u^2/2} \e^{-(x-\beta)^2/2}.
\end{multline}
Finally, using \eqref{Nexpo} and \eqref{OUcom} along with the strong Markov property, for any $R \ge R_2$ and any $x \ge \beta+1$,
\begin{align*}
&\sup_{\substack{z \in[-\beta, 0],\\y \ge 2R l_0(\beta)}}\prob_{(z,y)}\left(\tau_1(-x) < \tau_2(Rl_0(\beta))\right) = \sup_{\substack{z \in[-\beta, 0],\\y \ge 2R l_0(\beta)}}\prob_{(z,y)}\left(\inf_{t \le \sigma^{(\mathcal{N}^{\sigma})}}Q_1(t) < -x\right)\\
&\hspace{1cm}\le \sup_{\substack{z \in[-\beta, 0],\\y \ge 2R l_0(\beta)}}\sum_{k=0}^{\infty}\prob_{(z,y)}\left(\inf_{t \in [\sigma^{(2k+1)}, \sigma^{(2k+2)}]}Q_1(t) < -x, \mathcal{N}^{\sigma} \ge 2k+2\right)\\
&\hspace{1cm}\le \sup_{\substack{z \in[-\beta, 0],\\y \ge 2R l_0(\beta)}}\sum_{k=0}^{\infty}\mathbb{E}_{(z,y)}\mathbbm{1}_{[\mathcal{N}^{\sigma} \ge 2k+2]}\sup_{y > 0}\prob_{(-\beta-1, y)}\left(\tau_1(-x) < \tau_1(-\beta)\right)\\
&\hspace{1cm}\le \sup_{\substack{z \in[-\beta, 0],\\y \ge 2R l_0(\beta)}}\mathbb{E}_{(z,y)}(\mathcal{N}^{\sigma}) \ \sup_{y > 0}\prob_{(-\beta-1, y)}\left(\tau_1(-x) < \tau_1(-\beta)\right) \le C_2(\beta,R)\e^{-(x-\beta)^2/2}
\end{align*}
where $C_2(\beta,R)>0$ is a constant depending on $\beta, R$. This proves the lemma.
\end{proof}
\begin{lemma}\label{alphatwo}
For any $R > 1$ and any $x \ge 18Rl_0(\beta)$, there exists a constant $C_3(\beta,R)>0$ (depending on $\beta,R$) such that
\begin{equation*}
\sup_{z \in [-9Rl_0(\beta),0], \ y \le 2Rl_0(\beta)} \prob_{(z,y)}\left(\tau_1(-x) < \tau_2(2Rl_0(\beta))\right) \le C_3(\beta,R)\e^{-(x-\beta)^2/2}.
\end{equation*}
\end{lemma}
\begin{proof}
Fix any $R>1$, $Q_1(0) = z \ge -9Rl_0(\beta)$ and $Q_2(0)=y \le 2Rl_0(\beta)$. Define the stopping times: $\gamma^{(0)} =0$ and for $k \ge 0$,
\begin{align*}
\gamma^{(2k+1)} &= \inf\{t \ge \gamma^{(2k)}: Q_1(t) = -18Rl_0(\beta) \text{ or } Q_2(t) = 2Rl_0(\beta)\},\\
\gamma^{(2k+2)} &= \inf\{t \ge \gamma^{(2k+1)}: Q_1(t) = -9Rl_0(\beta) \text{ or } Q_2(t) = 2Rl_0(\beta)\}.
\end{align*}
Define $\mathcal{N}^{\gamma} = \inf\{n \ge 1: Q_2(\gamma^{(n)}) = 2Rl_0(\beta)\}$. Taking $B=2Rl_0(\beta)$ and $M=18Rl_0(\beta)$ in Lemma~\ref{lem:exit-by-q2}, we know there exists $q(\beta,R)$ such that
\begin{multline}\label{Nexpo2}
\inf_{z \in [-9Rl_0(\beta),0], \ y \le 2Rl_0(\beta)} \prob_{(z,y)}\left(\tau_2(2Rl_0(\beta)) < \tau_1(-18Rl_0(\beta))\right)\\
\ge \inf_{z \in [-9Rl_0(\beta),0], \ y \le 2Rl_0(\beta)} \prob_{(z,y)}\left(\tau_2(4Rl_0(\beta)) < \tau_1(-18Rl_0(\beta))\right) \ge q(\beta,R)>0.
\end{multline}
Using \eqref{Nexpo2} and the strong Markov property, there exists a constant $C(\beta,R)>0$ depending on $\beta,R$ such that
\begin{multline}\label{Nexpo3}
\sup_{z \in [-9Rl_0(\beta),0], \ y \le 2Rl_0(\beta)} \mathbb{E}_{(z,y)}(\mathcal{N}^{\gamma}) \le 2\sum_{n=0}^{\infty}\prob(\mathcal{N}^{\gamma} > 2n)\\
\le 2\sum_{n=0}^{\infty}(1-q(\beta,R))^n \le C(\beta,R) <\infty.
\end{multline}
Using \eqref{Nexpo3} and \eqref{OUcom} along with the strong Markov property, we obtain for any $x \ge 18Rl_0(\beta)$, 
\begin{align*}
&\sup_{\substack{z \in [-9Rl_0(\beta),0], \\ y \le 2Rl_0(\beta)}} \prob_{(z,y)}\left(\tau_1(-x) < \tau_2(2Rl_0(\beta))\right) \\
&\hspace{1cm}= \sup_{\substack{z \in [-9Rl_0(\beta),0], \\ y \le 2Rl_0(\beta)}} \prob_{(z,y)}\left(\inf_{t \le \gamma^{(\mathcal{N}^{\gamma})}}Q_1(t) < -x\right)\\
&\hspace{1cm}\le \sup_{\substack{z \in [-9Rl_0(\beta),0], \\ y \le 2Rl_0(\beta)}} \sum_{k=0}^{\infty}\prob_{(z,y)}\left(\inf_{t \in [\gamma^{(2k+1)}, \gamma^{(2k+2)}]}Q_1(t) < -x, \mathcal{N}^{\gamma} \ge 2k+2\right)\\
&\hspace{1cm}\le \sup_{\substack{z \in [-9Rl_0(\beta),0], \\ y \le 2Rl_0(\beta)}} \sum_{k=0}^{\infty}\mathbb{E}_{(z,y)}\mathbbm{1}_{[\mathcal{N}^{\gamma} \ge 2k+2]}\sup_{y > 0}\prob_{(-18Rl_0(\beta), y)}\left(\tau_1(-x) < \tau_1(-\beta)\right)\\
&\hspace{1cm}\le \sup_{\substack{z \in [-9Rl_0(\beta),0], \\ y \le 2Rl_0(\beta)}} \mathbb{E}_{(z,y)}(\mathcal{N}^{\gamma}) \ \sup_{y > 0}\prob_{(-18Rl_0(\beta), y)}\left(\tau_1(-x) < \tau_1(-\beta)\right)\\
&\hspace{1cm}\le C_3(\beta,R)\e^{-(x-\beta)^2/2}
\end{align*}
for some constant $C_3(\beta,R)>0$ depending on $\beta,R$. This proves the lemma.
\end{proof}
Now, we are in a position to give an upper bound to the fluctuations of $Q_1$ between two successive regeneration times $\Xi_k$ and $\Xi_{k+1}$, $k \ge 0$, defined in~\eqref{rendef} taking $B = Rl_0(\beta)$ for sufficiently large fixed $R$.
\begin{lemma}\label{Q1ub}
Fix any $R \ge \max\{2,R_1,R_2\}$, where $R_1$ and $R_2$ are obtained from Lemmas \ref{unifprob} and \ref{alphaone} respectively. 
Let $(Q_1(0), Q_2(0))=(0,2Rl_0(\beta))$ and take $B=Rl_0(\beta)$ in~\eqref{rendef}. 
There exists a constant $C^*(\beta,R)>0$ depending on $\beta,R$ such that for any $x \ge 18Rl_0(\beta)$,
\begin{equation*}
\prob_{(0,2Rl_0(\beta))}\left(\inf_{t \le \Xi_0} Q_1(t) < -x\right) \le C^*(\beta,R) \e^{-(x-2\beta)^2/8}.
\end{equation*}
\end{lemma} 
\begin{proof}
Choose and fix $R \ge \max\{2,R_1,R_2\}$. Define
$$
\Xi^* = \inf\big\{t \ge \tau_2(Rl_0(\beta)): Q_1(t) \ge -\beta -1\big\}.
$$
Then for any $x \ge 2(\beta+1)$,  by Lemma \ref{alphaone} and \eqref{OUcom} along with the strong Markov property,
\begin{equation}\label{interest}
\begin{split}
&\prob_{(0,2Rl_0(\beta))}\left(\inf_{\tau_2(Rl_0(\beta)) \le t \le \Xi^*} Q_1(t) < -x\right)\\
& \leq \prob_{(0,2Rl_0(\beta))}\left(\inf_{\tau_2(Rl_0(\beta)) \le t \le \Xi^*} Q_1(t) < -x, Q_1(\tau_2(Rl_0(\beta)))  \ge -x/2\right)\\
&\hspace{7.9cm}+ \prob_{(0,2Rl_0(\beta))}\left(\tau_1(-x/2) < \tau_2(Rl_0(\beta))\right)\\
&\le  \sup_{u \in [1, \frac{x}{2} - \beta]}\prob_{(-\beta-u, Rl_0(\beta))}\left(\tau_1(-x) < \tau_1(-\beta)\right)+\prob_{(0,2Rl_0(\beta))}\left(\tau_1(-x/2) < \tau_2(Rl_0(\beta))\right)\\
&\le  \sqrt{9\pi/2}\e^{(\frac{x}{2} - \beta)^2/2} \e^{-(x-\beta)^2/2}+C_2(\beta,R)\e^{-(\frac{x}{2}-\beta)^2/2}\\
&\le ( \sqrt{9\pi/2}+C_2(\beta,R))\e^{-(x-2\beta)^2/8}.
\end{split}
\end{equation}
Therefore, for any $x \ge 18Rl_0(\beta)$, using \eqref{interest} along with Lemmas~\ref{alphaone} and~\ref{alphatwo},
\begin{align*}
&\prob_{(0,2Rl_0(\beta))}\left(\inf_{t \le \Xi_0} Q_1(t) < -x\right) \le \prob_{(0,2Rl_0(\beta))}\left(\inf_{t \le \tau_2(Rl_0(\beta))} Q_1(t) < -x\right)\\
&\hspace{2cm} + \prob_{(0,2Rl_0(\beta))}\left(\inf_{\tau_2(Rl_0(\beta)) \le t \le \Xi^*} Q_1(t) < -x\right) + \prob_{(0,2Rl_0(\beta))}\left(\inf_{\Xi^* \le t \le \Xi_0} Q_1(t) < -x\right)\\
& \le \prob_{(0,2Rl_0(\beta))}\left(\tau_1(-x) < \tau_2(Rl_0(\beta))\right) + \prob_{(0,2Rl_0(\beta))}\left(\inf_{\tau_2(Rl_0(\beta)) \le t \le \Xi^*} Q_1(t) < -x\right)\\
&\hspace{2cm} + \sup_{z \in [-9Rl_0(\beta),0],\ y \le 2Rl_0(\beta)} \prob_{(z,y)}\left(\tau_1(-x) < \tau_2(2Rl_0(\beta))\right)\\
& \le C_2(\beta,R)\e^{-(x-\beta)^2/2} + ( C_2(\beta,R) + \sqrt{9\pi/2})\e^{-(x-2\beta)^2/8} + C_3(\beta,R)\e^{-(x-\beta)^2/2}\\
& \le C^*(\beta,R) \e^{-(x-2\beta)^2/8}
\end{align*}
which proves the lemma.
\end{proof}
Now, we prove a lower bound for the fluctuation of $Q_1$.
\begin{lemma}\label{lowqone}
Let $(Q_1(0), Q_2(0))=(0,2Rl_0(\beta))$ and and take $B=Rl_0(\beta)$ in~\eqref{rendef}. There exist constants $R^{**} >0$ not depending on $\beta$ such that for any $R \ge R^{**}$ and any $x \ge \beta$,
\begin{equation*}
\prob_{(0,2Rl_0(\beta))}\left(\inf_{t \le \Xi_0} Q_1(t) < -x\right) \ge C^{**}(\beta, R) \e^{-x^2},
\end{equation*}
where the positive constant $C^{**}(\beta, R)$ depends on both $\beta$ and $R$.
\end{lemma} 
\begin{proof}
Using $y = 2Rl_0(\beta) -\beta$ in Lemma \ref{lem:Q2hit1}, we observe that there exists $R^{**}>0$ such that for all $R \ge R^{**}$, there is a constant $q_1(\beta,R)>0$ (depending on $\beta,R$) for which
\begin{equation}\label{low1}
\prob_{(0,2Rl_0(\beta))}\left(\tau_1(-\beta/2) > \tau_2(Rl_0(\beta) + \beta/2\right) \ge q_1(\beta,R) >0.
\end{equation}
Recall $S(t) = Q_1(t) + Q_2(t)$. 
Recall that $Q_1(t) \le S(t) \le Q_2(t)$ for every $t$, and
when $Q_1(0)\in [0,\beta/2]$,
$$
S(t) = S(0) + \sqrt{2}W(t) - \beta t + \int_0^t(-Q_1(s))ds \le S(0) + \sqrt{2}W(t) - \frac{\beta}{2} t
$$
for $t \le \tau_1(-\beta/2)$.
Moreover, observe that if $Q_2(0)\leq 2Rl_0(\beta)$, then $Q_1(\tau_2(2Rl_0(\beta))) =0$ and consequently, $S(\tau_2(2Rl_0(\beta))) = Q_2(\tau_2(2Rl_0(\beta))) = 2Rl_0(\beta)$. 
Thus,
\begin{equation}\label{low2}
\begin{split}
&\sup_{z \in [-\beta/2,0]} \prob_{(z,Rl_0(\beta) + \beta/2)}\left(\tau_2(2Rl_0(\beta)) < \tau_1(-\beta/2) \right)\\
& \le \sup_{z \in [-\beta/2,0]} \prob_{(z,Rl_0(\beta) + \beta/2)}\left(S(t) \text{ hits } 2Rl_0(\beta) \text{ before } -\beta/2\right)\\
& \le \sup_{z \in [-\beta/2,0]} \prob_{(z,Rl_0(\beta) + \beta/2)}\left(z + Rl_0(\beta) + \beta/2 + \sqrt{2}W(t) -\frac{\beta}{2}t \text{ hits } 2Rl_0(\beta) \text{ before } -\beta/2\right)\\
& \le \prob\left(\sqrt{2}W(t) -\frac{\beta}{2}t \text{ hits } Rl_0(\beta) -\beta/2\right) \le \e^{-\beta(Rl_0(\beta) -\beta/2)/2}=: 1-q_2(\beta, R) <1.
\end{split}
\end{equation}
For $y \le 2Rl_0(\beta)$ and $(Q_1(0), Q_2(0)) = (-\beta/2, y)$, by \cite[Proposition 2.18]{Karatzas}, a process $U$ can be constructed on the same probability space as $ (Q_1, Q_2)$ such that almost surely $Q_1(t) +\beta \le U(t)$ for all $t \le \tau_1(0)$, where $U$ is an Ornstein-Uhlenbeck process which solves the SDE:
$$
dU(t) = \sqrt{2}dW(t) + (2Rl_0(\beta)-U(t))dt, \ \ U(0)=\beta/2.
$$
The scale function for $U$ is given by $s_U(u) = \int_0^u\e^{(w-2Rl_0(\beta))^2/2}dw$. Therefore, by elementary estimates on $s_U$, there exists a constant $C(\beta,R)>0$ (depending on $\beta,R$) such that for any $x \ge \beta$,
\begin{multline}\label{low3}
\inf_{y \le 2Rl_0(\beta)}\prob_{(-\beta/2, y)}\left(\tau_1(-x) < \tau_1(0)\right) \ge \prob\left(U(t) \text{ hits } -(x-\beta) \text{ before } \beta\right)\\
= \frac{s_U(\beta) - s_U(\beta/2)}{s_U(\beta) - s_U(-(x-\beta))} \ge C(\beta,R)\e^{-x^2}.
\end{multline}
Recall the notation $\sigma(t) = \inf\{s \ge t: Q_1(s)=0\}$ and define the stopping time 
$$\sigma_R = \inf\{t > \tau_2(Rl_0(\beta)+\beta/2): Q_2(t)=2Rl_0(\beta)\}.$$ 
From \eqref{low1}, \eqref{low2} and \eqref{low3} and the strong Markov property, for any $R \ge R^{**}$ and any $x \ge \beta$,
\begin{align*}
&\prob_{(0,2Rl_0(\beta))}\left(\inf_{t \le \Xi_0} Q_1(t) < -x\right)\\
&\ge \prob_{(0,2Rl_0(\beta))}\left(\tau_2(Rl_0(\beta)+\beta/2) < \tau_1(-\beta/2) < \sigma_R, \ \tau_1(-x) \in (\tau_1(-\beta/2), \sigma(\tau_1(-\beta/2)))\right)\\
&\ge \prob_{(0,2Rl_0(\beta))}\left(\tau_1(-\beta/2) > \tau_2(Rl_0(\beta) + \beta/2)\right) \times \inf_{z \in [-\beta/2,0]} \prob_{(z,Rl_0(\beta) + \beta/2)}\left(\tau_1(-\beta/2) < \tau_2(2Rl_0(\beta)) \right)\\
&\hspace{10cm}\times \inf_{y \le 2Rl_0(\beta)}\prob_{(-\beta/2, y)}\left(\tau_1(-x) < \tau_1(0)\right)\\
& \ge q_1(\beta,R)q_2(\beta, R)C(\beta,R)\e^{-x^2}.
\end{align*}
This proves the lemma.
\end{proof}

\begin{proof}[Proof of Theorem~\ref{th:excren}]
Fix any 
\begin{equation}\label{Rzerodef}
R_0 \ge 4\max\{64,\log(4\widetilde{C}_1)/\widetilde{C}_2, R_1,R_2, R^{**}\},
\end{equation}
where where $R_1, R_2$ and $R^{**}$ are obtained from Lemmas \ref{unifprob}, \ref{alphaone} and \ref{lowqone} respectively and $\widetilde{C}_1, \widetilde{C}_2$ are the constants defined in the statement of Lemma \ref{Q2gebeta1}. Choose $B = R_0l_0(\beta)$ in \eqref{rendef}.\\

\noindent
To prove (i), note that $y_0$ defined in Lemma \ref{Q2gebeta2} satisfies $y_0 + \beta < R_0l_0(\beta)$ for our specific choice of $R_0$. Therefore, taking $z=\frac{y-\beta}{2}$ in place of $y$ in  Lemma \ref{Q2gebeta2} and applying the strong Markov property at $\tau_2(z + \beta)$, we have for any $y \ge 4R_0l_0(\beta)$,
\begin{multline*}
\prob_{(0, 2R_0l_0(\beta))}\left(\tau_2(y) \le \Xi_0 \right) = \prob_{(0, 2R_0l_0(\beta))}\left(\tau_2(y) \le \tau_2(R_0l_0(\beta)) \right)\\
 \le \prob_{(0, z+\beta)}\left(\tau_2(2z+\beta) \le \tau_2(R_0l_0(\beta)) \right) \le \prob_{(0, z+\beta)}\left(\tau_2(2z+\beta) \le \tau_2(y_0 + \beta) \right)\le C^*_1 \e^{-C^*_2 \beta z}.
\end{multline*}
Part (ii) follows from Lemma \ref{Q2lb} by taking $B= R_0l_0(\beta)$. Parts (iii) and (iv) are direct consequences of Lemmas \ref{Q1ub} and \ref{lowqone} respectively.
\end{proof}

\section{Proofs of the main results}\label{sec:proofmain}

\begin{proof}[Proof of Theorem~\ref{th:statail}]
We will show that the tail bounds stated in the theorem hold with $C_R(\beta) = 18R_0l_0(\beta)$ and $D_R(\beta) = 4R_0l_0(\beta)$, where $R_0$ is defined in \eqref{Rzerodef} and $l_0(\beta)$ was defined in \eqref{lzerodef}. Taking $B= R_0l_0(\beta)$ in Theorem \ref{th:stationary}, note that for any $x\ge 0, y>0$,
\begin{align}\label{re}
\pi(Q_1(\infty) < -x) &= \frac{\mathbb{E}_{(0, 2R_0l_0(\beta))}\left(\int_{0}^{\Xi_0}\mathbbm{1}_{[Q_1(s) < -x]}ds\right)}{\mathbb{E}_{(0, 2R_0l_0(\beta))}\left(\Xi_0\right)},\nonumber\\
\pi(Q_2(\infty) > y)& = \frac{\mathbb{E}_{(0, 2R_0l_0(\beta))}\left(\int_{0}^{\Xi_0}\mathbbm{1}_{[Q_2(s) > y]}ds\right)}{\mathbb{E}_{(0, 2R_0l_0(\beta))}\left(\Xi_0\right)}.
\end{align}
To prove the theorem, we only need to estimate the numerators in the above representation. By the Cauchy-Schwarz inequality, for $x \ge 18R_0l_0(\beta)$,
\begin{multline*}
\mathbb{E}_{(0, 2R_0l_0(\beta))}\left(\int_{0}^{\Xi_0}\mathbbm{1}_{[Q_1(s) < -x]}ds\right) \le \mathbb{E}_{(0, 2R_0l_0(\beta))}\left(\mathbbm{1}_{[\tau_1(-x)]< \Xi_0]} (\Xi_0 - \tau_1(-x))\right)\\
 \le \sqrt{\prob_{(0,2R_0l_0(\beta))}(\tau_1(-x)< \Xi_0)} \sqrt{\mathbb{E}_{(0, 2R_0l_0(\beta))}(\Xi_0)^2} \le \sqrt{C^*(\beta)} \e^{-(x-2\beta)^2/16}\sqrt{\mathbb{E}_{(0, 2R_0l_0(\beta))}(\Xi_0^2)},
\end{multline*}where the last inequality is a consequence of Part (iii) of Theorem~\ref{th:excren}. By Proposition \ref{prop:Xi}, $\mathbb{E}_{(0,2R_0l_0(\beta))}\left(\Xi^2_0\right) < \infty$. 
Now, using this in the above bound, we obtain the upper bound on $\pi(Q_1(\infty) < -x)$ claimed in the theorem. The upper bound for $\pi(Q_2(\infty)>y)$ is obtained similarly using part (i) of Theorem \ref{th:excren}.

To obtain the lower bound on $\pi(Q_1(\infty) < -x)$, we proceed along the same line of arguments as in the proof of Lemma \ref{lowqone}. Recall the stopping time 
$$\sigma_R = \inf\{t > \tau_2(Rl_0(\beta)+\beta/2): Q_2(t)=2Rl_0(\beta)\}.$$ Observe that for $x \ge \beta$,
\begin{align*}
&\mathbb{E}_{(0, 2R_0l_0(\beta))}\left(\int_{0}^{\Xi_0}\mathbbm{1}_{[Q_1(s) < -x]}ds\right)\\
& \ge \mathbb{E}_{(0, 2R_0l_0(\beta))}\left(\mathbbm{1}_{[\tau_2(R_0l_0(\beta) + \beta/2) < \tau_1(-\beta/2) < \sigma_{R_0}]}\int_{\tau_1(-\beta/2)}^{\sigma(\tau_1(-\beta/2))}\mathbbm{1}_{[Q_1(s) < - x]}ds\right)\\
& \ge \prob_{(0, 2R_0l_0(\beta))}\left(\tau_2(R_0l_0(\beta) + \beta/2) < \tau_1(-\beta/2)\right) \times\inf_{y \le 2R_0l_0(\beta)}\mathbb{E}_{(-\beta/2,y)}\left(\int_0^{\tau_1(0)}\mathbbm{1}_{[Q_1(s) < - x]}ds\right) \\
&\hspace{6cm} \times \inf_{z \in [-\beta/2,0]} \prob_{(z,R_0l_0(\beta) + \beta/2)}\left(\tau_1(-\beta/2) < \tau_2(2R_0l_0(\beta)) \right)\\
& \ge q_1(\beta,R_0)q_2(\beta, R_0)\inf_{y \le 2R_0l_0(\beta)}\mathbb{E}_{(-\beta/2,y)}\left(\int_0^{\tau_1(0)}\mathbbm{1}_{[Q_1(s) < - x]}ds\right)
\end{align*}
where $q_1(\beta,R_0)>0, q_2(\beta,R_0)>0$ are obtained in \eqref{low1} and \eqref{low2} respectively with $R_0$ in place of $R$.

Recall that for $y \le 2R_0l_0(\beta)$ and $(Q_1(0), Q_2(0)) = (-\beta/2, y)$, by  \cite[Proposition 2.18]{Karatzas}, a process $U_{\beta/2}$ can be constructed on the same probability space as $(Q_1,Q_2)$, such that $Q_1(t) +\beta \le U_{\beta/2}(t)$ for $t \le \tau_1(0)$, where $U_z$ is an Ornstein-Uhlenbeck process which solves the SDE:
$$
dU_z(t) = \sqrt{2}dW(t) + (2Rl_0(\beta)-U(t))dt, \ \ U_z(0)=z.
$$
where the scale function for $U_z$ is given by $s_U(u) = \int_0^u\e^{(w-2Rl_0(\beta))^2/2}dw$.

Define $\tau_z^U(w) = \inf\{t \ge 0: U_z(t)=w\}$ and write the law of $U_z$ and the corresponding expectation as $\prob^U_z$ and $\mathbb{E}^U_z$ respectively. Then, for $x \ge \beta$,
\begin{equation}\label{lowstatone}
\begin{split}
&\inf_{y \le 2R_0l_0(\beta)}\mathbb{E}_{(-\beta/2,y)}\left(\int_0^{\tau_1(0)}\mathbbm{1}_{[Q_1(s) < - x]}ds\right) \ge \mathbb{E}^U_{\beta/2}\left(\int_0^{\tau_{\beta/2}^U(\beta)}\mathbbm{1}_{[U_{\beta/2}(s) < - x+\beta]}ds\right)\\
&\ge \prob^U_{\beta/2}\left(\tau_{\beta/2}^U(-2x + \beta) < \tau_{\beta/2}^U(\beta)\right)\mathbb{E}^U_{-2x+\beta}\left(\tau_{-2x+\beta}^{U}(-x+\beta)\right), \quad\text{ by strong Markov property},\\
&= \frac{s_U(\beta) - s_U(\beta/2)}{s_U(\beta) - s_U(-(2x-\beta))} \ \mathbb{E}^U_{-2x+\beta}\left(\tau_{-2x+\beta}^{U}(-x+\beta)\right)
 \ge C(\beta) \e^{-4x^2}\mathbb{E}^U_{-2x+\beta}\left(\tau_{-2x+\beta}^{U}(-x+\beta)\right)
\end{split}
\end{equation}
where $C(\beta)$ is a positive constant that only depends on $\beta$. 
Now, from the Doob representation of Ornstein-Uhlenbeck process,
$$
U_{-2x+\beta}(t) = (-2x+\beta)\e^{-t} + 2R_0l_0(\beta)(1-\e^{-t}) + \e^{-t}\widetilde{W}(\e^{2t}-1)
$$
for a standard Brownian motion $\widetilde{W}$.
Therefore, taking $T=\log(5/4)$, for $x \ge 4R_0 l_0(\beta)$,
\begin{align*}
&\prob^U_{-2x+\beta}\left(\tau_{-2x+\beta}^U(-x + \beta)  \le T\right)\\
& \le \prob\left(\sup_{t \le T}\left((-2x+\beta)\e^{-t} + 2R_0l_0(\beta)(1-\e^{-t}) + \e^{-t}\widetilde{W}(\e^{2t}-1)\right) > -x+\beta\right)\\
& \le \prob\left((-2x+\beta)\e^{-T} + 2R_0l_0(\beta)(1-\e^{-T}) + \sup_{t \le T}\left(\widetilde{W}(\e^{2t}-1)\right) > -x+\beta\right)\\
& \le \prob\left(\sup_{t \le T}\left(\widetilde{W}(\e^{2t}-1)\right) > x/2\right) \ \ \text{ by our choice of $T$}\\
& = \prob\left(\sup_{t \le 1}\widetilde{W}(t) > \frac{x}{2\sqrt{\exp(2T)-1}}\right) \ \text{ by Brownian scaling }
 \le \frac{4\sqrt{\exp(2T)-1}}{\sqrt{2\pi}x} < \frac{1}{2}.
\end{align*}
Thus,
$$
\mathbb{E}^U_{-2x+\beta}\left(\tau_{-2x+\beta}^{U}(-x+\beta)\right) =\int_0^{\infty}\prob^U_{-2x+\beta}\left(\tau_{-2x+\beta}^U(-x + \beta) > t\right)dt \ge \frac{1}{2}\log(5/4).
$$
Using this in \eqref{lowstatone} gives us the lower bound on $\pi(Q_1(\infty) < -x)$ claimed in the theorem.

Finally, we prove the lower bound on $\pi(Q_2(\infty) > y)$.
Note that by the strong Markov property, for any $y \ge R_0l_0(\beta)$,
\begin{equation}\label{lowtwo}
\begin{split}
\mathbb{E}_{(0, 2R_0l_0(\beta))}\left(\int_{0}^{\Xi_0}\mathbbm{1}_{[Q_2(s) > y]}ds\right) &\ge \prob_{(0, 2R_0l_0(\beta))}\left(\tau_2(2y) \le \Xi_0\right) \times \mathbb{E}_{(0,2y)}\left(\tau_2(y)\right)\\
 &\ge (1-\e^{-\beta R_0l_0(\beta)})\e^{-\beta(2y-2R_0l_0(\beta))}\mathbb{E}_{(0,2y)}\left(\tau_2(y)\right)
\end{split}
\end{equation}
where the last step follows from Part (ii) of Theorem~\ref{th:excren}. Recall that
$$
Q_2(t) \ge S(t) \ge S(0) + \sqrt{2}W(t) - \beta t, \ \ t \ge 0,
$$
where $S(t) = Q_1(t) + Q_2(t)$. Therefore, stsrting with $(Q_1(0),Q_2(0)) = (0,2y)$, the hitting time of level $y$ of $Q_2$ is stochastically bounded below by the hitting time of $y$ by $S(0) + \sqrt{2}W(t) - \beta t$. Denoting the latter hitting time by $\tau^S(y)$, we obtain $\mathbb{E}_{(0,2y)}\left(\tau_2(y)\right) \ge \mathbb{E}_{(0,2y)}\left(\tau^S(y)\right)$. For $y \ge R_0 l_0(\beta)$,
\begin{multline*}
\prob_{(0,2y)}\left(\tau^S(y) \le \frac{y}{2\beta}\right) = \prob\left(\inf_{t \le \frac{y}{2\beta}}\left(2y + \sqrt{2}W(t) - \beta t\right) < y\right) \le \prob\left(\inf_{t \le \frac{y}{2\beta}}\left(\sqrt{2}W(t)\right) < -y/2\right)\\
= \prob\left(\inf_{t \le 1}\left(W(t)\right) < -\sqrt{\beta y}/2\right) \le \frac{4}{\sqrt{2\pi \beta y}} \le \frac{4}{\sqrt{2\pi R_0}} < \frac{1}{2},
\end{multline*}
for our choice of $R_0$. This gives
$$
\mathbb{E}_{(0,2y)}\left(\tau^S(y)\right) = \int_0^{\infty} \prob_{(0,2y)}\left(\tau^S(y) > t\right) dt \ge \frac{y}{4\beta}.
$$
Using this in \eqref{lowtwo} gives us the lower bound on $\pi(Q_2(\infty) > y)$ claimed in the theorem.
\end{proof}  

\begin{proof}[Proof of Theorem~\ref{th:lil}]
Below we provide the proof of the fluctuation result for $Q_2$.
The proof for $Q_1$ follows using analogous arguments.

Take $\mathcal{C^*}$ in the theorem to be the positive constant $C^*_2$ not depending on $\beta$ that was obtained in Part (i) of Theorem~\ref{th:excren}. 
Fix $\epsilon \in (0,1/2)$. Fix any starting point $(Q_1(0), Q_2(0))=(x,y)$. Then by Parts (i) and (ii) of Theorem~\ref{th:excren}, we obtain constants $D_1(\beta)$ and $D_2(\beta)$ and an integer $N(\beta)>0$ depending only on $\beta$ and an such that for all $n \ge N(\beta)$,
\begin{eqnarray*}
\prob_{(x,y)}\left(\sup_{t \in [\Xi_{n}, \Xi_{n+1}]} Q_2(t) > \frac{2(1+\epsilon)\log n}{C^*_2 \beta}\right) &\le \frac{D_1(\beta)}{n^{1+\epsilon}},\\
\prob_{(x,y)}\left(\sup_{t \in [\Xi_{n}, \Xi_{n+1}]} Q_2(t) > \frac{(1-\epsilon)\log n}{\beta}\right) &\ge \frac{D_2(\beta)}{n^{1-\epsilon}}.
\end{eqnarray*}
Therefore, by the Borel-Cantelli Lemma,
\begin{equation}\label{BC}
\frac{1-\epsilon}{\beta} \le \limsup_{n \rightarrow \infty}\frac{\sup_{t \in [\Xi_{n}, \Xi_{n+1}]} Q_2(t)}{\log n} \le \frac{2(1+\epsilon)}{C^*_2 \beta}, \ \ a.s.
\end{equation}
By Proposition \ref{prop:Xi}, $\mathbb{E}_{(0,2R_0l_0(\beta))}\left(\Xi_0\right) < \infty$ and as $\{\Xi_{n+1} - \Xi_{n}\}_{n \ge 0}$ are i.i.d., therefore by the Strong Law of Large Numbers,
\begin{equation}\label{SLLN}
\lim_{n \rightarrow \infty} \frac{\Xi_n}{n} \rightarrow \mathbb{E}_{(0,2R_0l_0(\beta))}\left(\Xi_0\right), \ \ a.s.
\end{equation}
From the lower bound in \eqref{BC}, with probability one, there exists a subsequence $\{n_k\}\subseteq\{n\}$ and $t_{n_k} \in [\Xi_{n_k}, \Xi_{n_k+1}]$ such that 
$$
Q_2(t_{n_k}) \ge (1-2\epsilon)\frac{\log n_k}{\beta}
$$ 
for all sufficiently large $k$. Moreover, by \eqref{SLLN}, almost surely,
$$
\log t_{n_k} \le \log \Xi_{n_k + 1} = \log \left(\frac{\Xi_{n_k+1}}{n_k+1}\right) + \log (n_k+1) \le (1+\epsilon)\log n_k
$$
for all sufficiently large $k$. Therefore, almost surely, for all sufficiently large $k$,
$$
\frac{Q_2(t_{n_k})}{\log t_{n_k}} \ge \frac{1-2\epsilon}{(1+\epsilon)\beta}.
$$
Since this holds for every $\epsilon \in (0,1/2)$, we obtain
$$
\limsup_{t \rightarrow \infty} \frac{Q_2(t)}{\log t} \ge \frac{1}{\beta}, \ \ a.s.
$$
From the upper bound in \eqref{BC} and \eqref{SLLN}, we obtain $n_0$ such that for all $n \ge n_0$ 
$$
\frac{\sup_{t \in [\Xi_{n}, \Xi_{n+1}]} Q_2(t)}{\log n} \le \frac{2(1+\epsilon)}{C^*_2 \beta}, \ \ \text{ and }\ \log t \ge (1-\epsilon) \log n.
$$
Therefore,
$$
\frac{Q_2(t)}{\log t} \le \frac{2(1+\epsilon)}{(1-\epsilon)C^*_2 \beta}, \ \ \text{ for all } t \ge \Xi_{n_0}
$$
and hence,
$$
\limsup_{t \rightarrow \infty} \frac{Q_2(t)}{\log t} \le \frac{2}{C^*_2 \beta}, \ \ a.s.
$$
The fluctuation result for $Q_1$ is obtained similarly using Parts (iii) and (iv) of Theorem \ref{th:excren}.
\end{proof}

\section*{Acknowledgement}
The authors sincerely thank Amarjit Budhiraja for several insightful discussions.
DM was supported by The Netherlands Organization for Scientific Research (NWO) through Gravitation Networks grant 024.002.003, and TOP-GO grant 613.001.012.
The work was initiated during DM's visit to UNC, Chapel hill. 
DM sincerely thanks the hospitality of UNC, Chapel hill for that.


\begin{appendices}
\section{Proof of Lemma \ref{lem:q2regeneration}}\label{app:lem4.3}

In this appendix we will prove Lemma~\ref{lem:q2regeneration}.
As mentioned earlier, we need to have sharp estimates for the time  $Q_2$ takes to hit the level $B$ starting from a large initial state.
This, in turn, amounts to estimating the time integral of the $Q_1$ process when $Q_2$ is large, 
which is furnished by Lemma~\ref{lem:q1integral}.
The tail estimates presented in Lemmas~\ref{lem:integrated} and~\ref{lem:integrated2} will be used in the proof of Lemma~\ref{lem:q1integral}.

Fix any $M>0$ and $\varepsilon>0$. 
Observe that if $\inf_{0\leq s\leq t}Q_2(s)>M + \beta$, then the process $\{Q_1(s)\}_{0\leq s\leq t}$ is bounded below by the process $\{\eta(s)\}_{0\leq s\leq t}$, where 
$$\eta(t) = Q_1(0) + \sqrt{2}W(t) + Mt - L_\eta(t),$$ 
with $L_\eta$ being
the local time of $\eta$ given by $L_{\eta}(t) = \sup_{s \le t}\{Q_1(0) + \sqrt{2}W(s) + Ms\}^+$ (where $x^+ = \max\{x,0\}$ for any $x \in \mathbb{R}$), and $W$ being the standard Brownian motion.
Note that the dependence of $M$ in $\eta$ is suppressed for convenience in notation.
 For $i\geq 1$ define 
\begin{align*}
&T_{2i-1}:= \inf\ \{t>T_{2i-2}: \eta(t) = -\varepsilon\},
&T_{2i}:= \inf\ \{t>T_{2i-1}: \eta(t) = -\varepsilon/2\},\\
&\xi_i:= T_{2i}-T_{2i-1},\quad	\zeta_i := T_{2i+1}-T_{2i},
&u_i:= \sup_{T_{2i-1}\leq t\leq T_{2i}}(-\eta(t)),\\
&N_t = \inf\ \{n\geq 1: T_{2n}\geq t\}.
\end{align*}
with the convention that $T_{0}\equiv 0.$
Further, for $i\geq 1$, let $T_i^W$ denote the corresponding stopping
times when the process $\eta$ is replaced by the process $W_R$ described as 
$$W_R(t) = Q_1(0) + \sqrt{2}W(t) - L_W(t)$$
with $L_W$ being
the local time of $W_R$ given by $L_W(t) = \sup_{s \le t}\{Q_1(0) + \sqrt{2}W(s)\}^+$.
Also, similarly denote $\xi_i^W:= T_{2i}^W-T_{2i-1}^W$ and $\zeta_i^W := T_{2i+1}^W-T_{2i}^W$. 
\begin{lemma}\label{lem:integrated}
Assume that $Q_1(0) \in [-\varepsilon, 0]$.
Then the following hold:
\begin{enumerate}[{\normalfont (i)}]
\item For $i\geq 1$, $\zeta_i^W\leq_{st}\zeta_i$.
\item \label{subexp} There exist constants $c_\zeta, c_W>0$ not depending on $M, \varepsilon$ such that for $t \ge \varepsilon^2$
\begin{equation}
\begin{split}
\text{(a)}~\Pro{\zeta_1>t}\geq \exp(-c_\zeta t/\varepsilon^2) \qquad\text{and}\qquad \text{(b)}~\Pro{\zeta_1^W>t}\leq \exp(-c_Wt/\varepsilon^2).
\end{split}
\end{equation}
\item\label{lem:tails}
For all $x\geq \varepsilon$,
$\Pro{u_1>x}\leq \exp(-M(x-\varepsilon)),$
\item For all $t\geq \varepsilon/M$,
$\Pro{\xi_1>t}\leq \frac{2}{\sqrt{\pi}M\sqrt{t}}\exp(-M^2 t/16)$.
\item \label{lem:Nt-upper}
There exist constants $b,c_N^{(1)}>0$ not depending on $M, \varepsilon$, such that for $t \ge \varepsilon^2/b$
$$\Pro{N_t > b\varepsilon^{-2}t}\leq 2\exp(-c_N^{(1)}t/\varepsilon^2).$$ 
\end{enumerate}
\end{lemma}
\begin{proof}
(i) This is an immediate consequence of the fact that $\{\eta(s)\}_{0\leq s\leq t}\geq_{st} \{W_R(s)\}_{0\leq s\leq t}$. \\

\noindent
(ii)
Take $\varepsilon=1$. Using the Markov property for reflected Brownian motion, it is easy to see that there exist constants $c_\zeta, c_W>0$ such that $\exp(-c_Wt)\geq\Pro{\zeta_1^W>t}\geq \exp(-c_\zeta t)$ for $t \ge 1$.
(ii.a) now follows from (i) and Brownian scaling.
(ii.b) is also an immediate consequence of Brownian scaling. \\

\noindent
(iii) Observe that
\begin{align*}
\Pro{u_1>x} 
\leq \Pro{\inf_{s<\infty}(-\varepsilon + \sqrt{2}W(s) + Ms)<-x} = \exp(-M(x-\varepsilon)),
\end{align*}
since $-\inf_{s < \infty}(\sqrt{2}W(s) + Ms)$ follows an exponential random variable with mean $1/M$.\\

\noindent
(iv) Note that
\begin{align*}
\Pro{\xi_1>t} &= \Pro{\sup_{s\leq t}(-\varepsilon+\sqrt{2}W(s) +Ms) \le -\varepsilon/2}\\
&\leq \Pro{\sqrt{2}W(t) + Mt \leq \varepsilon/2}
\leq \frac{2}{\sqrt{\pi}M\sqrt{t}}\exp(-M^2 t/16)\qquad\forall\ t\geq \varepsilon/M.
\end{align*}

\noindent
(v)
Observe that
\begin{align*}
\Pro{N_t>b\varepsilon^{-2}t} &\leq \Pro{\sum_{i=1}^{\lfloor b\varepsilon^{-2}t \rfloor}\zeta_i\leq t}\leq \Pro{\sum_{i=1}^{\lfloor b\varepsilon^{-2}t \rfloor}\zeta_i^W\leq t}, \ \text{ by part (i),}\\
&\leq \prob\Big(\sum_{i=1}^{\lfloor b\varepsilon^{-2}t \rfloor}\varepsilon^{-2}\left(\zeta_i^W - \expt(\zeta_i^W)\right)\leq - \left(\frac{b}{2}\varepsilon^{-2}\expt \zeta_1^W - 1\right)t\varepsilon^{-2}\Big)\\
&\leq 2\exp(-c_N^{(1)}t/\varepsilon^2) \qquad\mbox{[choosing }b = 4\varepsilon^2/\expt(\zeta_1^W)\mbox{]},
\end{align*}
where the last step follows from part (ii), which shows that $\varepsilon^{-2}\left(\zeta_i^W - \expt(\zeta_i^W)\right)$ are sub-exponential random variables, and then using the Chernoff's inequality (see \cite[Pg.~16, Equation (2.2)]{massart2007}) to the sum $\sum_{i=1}^{\lfloor b\varepsilon^{-2}t \rfloor}\varepsilon^{-2}\left(\zeta_i^W - \expt(\zeta_i^W)\right)$.
Here, note that by Brownian scaling, $b$ chosen above does not depend on $\varepsilon$.
\end{proof}

The next technical lemma establishes a useful concentration inequality that will be crucial in obtaining tail probabilities for $\sum_{i=1}^{N_t}u_i\xi_i$.
\begin{lemma}\label{lem:BM-concentration}
Fix $\varepsilon>0$ and $M \ge \frac{1}{\varepsilon}$. 
Let $\Phi_i$'s be iid nonnegative random variables with 
$$\Pro{\Phi_1>z}\leq \exp(-c'M^{3/2}\sqrt{z})\quad \text{for all}\quad z\geq 4\varepsilon^2/M,$$ 
and $\E{\Phi_1}\leq c_{11}\varepsilon^2/M$ where $c', c_{11}$ are positive constants not depending on $M, \varepsilon$. Then 
$$\prob\Big(\sum_{i=1}^n\Phi_i\geq 4c_{11}n\frac{\varepsilon^2}{M}\Big)\leq \Big(1 + c_1\frac{1}{n^{2/5}\left(\varepsilon M\right)^{8/5}}\Big)\exp\Big(-c_2 (\varepsilon M)^{4/5}n^{1/5}\Big),$$
for $n \ge c_3\varepsilon M$, where $c_1,c_2, c_3$ are positive constants not depending on $M, \varepsilon$.
\end{lemma}
\begin{proof}
For some $A \ge 4\varepsilon^2/M$ to be chosen later, define 
$$\Phi_i^*:= \Phi_i\ind{\Phi_i\geq A} \qquad\mbox{and}\qquad \Phi_i^{**}:= \Phi_i\ind{\Phi_i < A}.$$
Thus, $\Phi_i = \Phi_i^* + \Phi_i^{**}$. Note that
\begin{align*}
\E{\Phi_i^*}^2 &= \int_{A^2}^{\infty}\prob\left(\Phi_i>\sqrt{z}\right)dz = \int_{A}^{\infty}2z\prob\left(\Phi_i>z\right)dz \\
&\le  \int_{A}^{\infty}2z\exp(-c'M^{3/2}\sqrt{z})dz
= \int_{\sqrt{A}}^{\infty}4z^3\exp(-c'M^{3/2}z)dz\\
& =  \frac{4}{M^6}\int_{M^{3/2}\sqrt{A}}^{\infty}z^3\exp(-c'z)dz
\leq c'' \frac{A^{3/2}}{M^{3/2}} \exp(-c'M^{3/2}\sqrt{A}),
\end{align*}
where the constant $c''$ does not depend on $M,A$.
Thus, using Chebyshev's inequality, 
\begin{equation}\label{eq:*bound}
\Pro{\sum_{i=1}^n\Phi^*>2c_{11}n\frac{\varepsilon^2}{M}}\leq \frac{c'' M^{1/2} A^{3/2}\exp(-c'M^{3/2}\sqrt{A})}{4nc_{11}^2 \varepsilon^4}.
\end{equation}
Further note that $\Phi_i^{**}$'s are bounded random variables.
Therefore using Azuma-Hoeffding inequality we obtain, 
\begin{equation}\label{eq:**bound}
\begin{split}
\Pro{\sum_{i=1}^n\Phi_i^{**}>2c_{11}n\frac{\varepsilon^2}{M}}
&= \Pro{\sum_{i=1}^n(\Phi_i^{**} - \E{\Phi_i^{**}})>c_{11}n\frac{\varepsilon^2}{M}}\\
&\leq \exp\Big(-\big(\frac{c_{11}n\varepsilon^2}{M}\big)^2/(8A^2n)\Big)= \exp(-c_{11}^2n\varepsilon^4/(8A^2M^2))
\end{split}
\end{equation}
Equating the exponents of equations~\eqref{eq:*bound} and~\eqref{eq:**bound}, and solving for $A$, we get
$$A = \left(\frac{c_{11}^2}{8c'}\right)^{2/5}\left(\frac{\varepsilon^{8/5}n^{2/5}}{M^{7/5}}\right).$$ The condition $A \ge 4\varepsilon^2/M$ implies $n \ge 2^5\left(\frac{8c'}{c_{11}^2}\right)\varepsilon M$. This choice for $A$ yields the bound claimed in the lemma. 
\end{proof}

\begin{lemma}\label{lem:integrated2}
Fix any $\varepsilon>0$ and $M \ge \frac{1}{\varepsilon}$.
\begin{enumerate}[{\normalfont (i)}]
\item \label{lem:u1xi1} There exist positive constants $c', c_{11}$ not depending on $M, \varepsilon$, such that
\begin{align*}
\mbox{\normalfont (a)}&\quad \Pro{u_1\xi_1>x} \le \exp(-c'M^{3/2}\sqrt{x})\qquad \forall \ x \ge 4\varepsilon^2/M,\\
\mbox{\normalfont (b)}&\quad \E{u_1\xi_1}\leq c_{11}\frac{\varepsilon^2}{M}.
\end{align*}
\item \label{fact:sum-uixi}
Let $b,c_{11}$ be the constants in Lemma~\ref{lem:integrated} \eqref{lem:Nt-upper} and Lemma~\ref{lem:integrated2} \eqref{lem:u1xi1} respectively. 
There exist constants $c_1, c_2, c_3$ not depending on $\varepsilon, M$ such that
$$
\Pro{\sum_{i=1}^{N_t}u_i\xi_i> 4 \frac{bc_{11}t}{M}} \leq c_1\exp(-c_2(\varepsilon M)^{4/5}(t/\varepsilon^2)^{1/5})
$$
for $t \ge c_3 \varepsilon^3M$.
\end{enumerate}
\end{lemma}
\begin{proof}
(i.a)~Recall that $M \ge \frac{1}{\varepsilon}$.
By Lemma~\ref{lem:integrated}~\eqref{lem:tails}, we obtain for $x \ge 4\varepsilon^2/M,$
\begin{align*}
\Pro{u_1\xi_1>x}&\leq \Pro{u_1>\sqrt{Mx}}+\Pro{u_1\xi_1>x, u_1\leq \sqrt{Mx}}\\
&\leq \Pro{u_1>\sqrt{Mx}}+\Pro{\xi_1>\frac{\sqrt{x}}{\sqrt{M}}}\\
&\leq \exp(-M(\sqrt{Mx}-\varepsilon))+\frac{2}{\sqrt{\pi}M^{3/4}x^{1/4}}\exp(-M^{3/2} \sqrt{x}/16)\\
&\leq \exp(-M^{3/2}\sqrt{x}/2)+\frac{2}{\sqrt{\pi}M^{3/4}x^{1/4}}\exp(-M^{3/2} \sqrt{x}/16)\\
&\leq \exp(-c'M^{3/2}\sqrt{x}),
\end{align*}
where the last line is a consequence of the fact that for $x \ge 4\varepsilon^2/M$ and $M \ge \frac{1}{\varepsilon}$, $M^{3/4}x^{1/4} \ge \sqrt{2M\varepsilon} > 1$.\\

\noindent
(i.b)~As a consequence of part (i.a) we obtain
\begin{align*}
\E{u_1\xi_1}&\leq \int_0^{4\varepsilon^2/M}\dif x
+\frac{1}{M^3}\int_{4\varepsilon^2/M}^\infty \exp(-c'M^{3/2}\sqrt{x})M^3\dif x\\
&\leq \frac{4\varepsilon^2}{M} + \frac{c'''}{M^3} \leq c_{11}\frac{\varepsilon^2}{M},
\end{align*}
where we again used $M \ge \frac{1}{\varepsilon}$ to obtain $\frac{1}{M^3} \le \frac{\varepsilon^2}{M}$.\\

\noindent
(ii)
Observe that due to Lemma~\ref{lem:integrated} \eqref{lem:Nt-upper} and Lemma~\ref{lem:BM-concentration},
\begin{align*}
\prob\Big(\sum_{i=1}^{N_t}u_i\xi_i>4 \frac{bc_{11}t}{M}\Big)
&\leq\Pro{N_t>b\varepsilon^{-2}t} + \prob\ \Big(\sum_{i=1}^{\lfloor b\varepsilon^{-2}t\rfloor}u_i\xi_i>4 \frac{bc_{11}t}{M}\Big)\\
&\leq 2\exp(-c_N^{(1)}t/\varepsilon^2) + C_1\exp(-C_2(\varepsilon M)^{4/5}(t/\varepsilon^2)^{1/5})
\end{align*}
for $t \ge C_3 \varepsilon^3M$, where $C_1, C_2, C_3$ can be chosen to be independent of $M, \varepsilon$. This completes the proof.
\end{proof}

We are now in a position to state and prove Lemma~\ref{lem:q1integral} that provides us with a crucial estimate for the time-integral of the $Q_1$ process when $Q_2$ is large.
\begin{lemma}\label{lem:q1integral}
There exist $c'_1, c'_2, c'_3>0$, not depending on $\beta$ such that for any $y > c'_1\left(\beta \vee \beta^{-1}\right)+ \beta$,
\begin{align*}
&\prob_{(0, y)}\Bigg(\int_0^t(-Q_1(s))\dif s>\left(\beta \wedge \beta^{-1}\right)\frac{t}{2},\ \inf_{s\leq t}Q_2(s) \geq c'_1\left(\beta \vee \beta^{-1}\right)+ \beta\Bigg)\\
&\leq \exp\Big(-c'_2t^{1/5}\left(\beta \vee \beta^{-1}\right)^{2/5}\Big)\qquad \text{for}\quad t \ge c'_3 \left(\beta \wedge \beta^{-1}\right)^2.
\end{align*} 
\end{lemma}
\begin{proof}
Recall the constants $b$ and $c_{11}$ from Lemma~\ref{lem:integrated}~\eqref{lem:Nt-upper} and Lemma~\ref{lem:integrated2}~\eqref{lem:u1xi1} respectively. 
As $c_{11}$ appears in the upper bound of $\mathbb{E}\left(u_1\xi_1\right)$ in Lemma~\ref{lem:integrated2}~\eqref{lem:u1xi1}, we can take $c_{11}> b^{-1} \vee 1$. First we consider the case $\beta \in (0,1)$. Take $\varepsilon = \beta/4$. Choose $M = 16c_{11}b/\beta$, since in that case 
$$\varepsilon = \frac{\beta}{4} = \frac{4c_{11}b}{M}.$$
Observe that 
\begin{align*}
&\prob_{(0,y)}\Big(\int_0^t(-Q_1(s))\dif s>\frac{\beta t}{2},\ \inf_{s\leq t}Q_2(s) \geq M+\beta\Big)\\
&\leq\prob_{(0,\ y)}\Big(\sum_{i=1}^{N_t}\int_{T_{2i-1}}^{T_{2i}}(-Q_1(s))\dif s > \frac{4c_{11}b}{M}t, \inf_{s\leq t}Q_2(s) \geq M+\beta\Big)\\
&\leq\Pro{\sum_{i=1}^{N_t}u_i\xi_i>\frac{4c_{11}b}{M}t}\\
&\leq \exp\Big(-c''_2(\beta M)^{4/5}(t/\beta^2)^{1/5}\Big)\\
&\leq \exp\Big(-c'_2(t/\beta^2)^{1/5}\Big) \qquad \text{for}\quad t \ge c''_3 \beta^3 M= c'_3\beta^2, \qquad\mbox{due to Lemma~\ref{lem:integrated2}~\eqref{fact:sum-uixi}},
\end{align*}
where the constants $c'_2, c''_2 c'_3, c''_3$ do not depend on $\beta,M$. 
Next, for the case $\beta>1$, we take $\varepsilon = \frac{1}{4\beta}$ and $M= 16c_{11}\beta$ so that
$$
\varepsilon = \frac{1}{4\beta} = \frac{4c_{11}b}{M},
$$
and then apply the same argument. This completes the proof.
\end{proof}

\begin{proof}[Proof of Lemma~\ref{lem:q2regeneration}]
Let us denote the following events
\begin{align*}
\cE_t&:= \Big[\inf_{s\leq t} Q_2(s) > c'_1\left(\beta \vee \beta^{-1}\right)+\beta\Big],\\
\cE_t^1&:= \Bigg[\int_0^t(-Q_1(s))\dif s>\frac{\beta t}{2}, \quad\inf_{s\leq t} Q_2(s) > c'_1\left(\beta \vee \beta^{-1}\right)+ \beta\Bigg],\\
\cE_t^2&:=\Bigg[\int_0^t(-Q_1(s))\dif s\leq\frac{\beta t}{2},\quad \inf_{s\leq t} Q_2(s) > c'_1\left(\beta \vee \beta^{-1}\right)+ \beta\Bigg].
\end{align*}
Note that if $(Q_1(0), Q_2(0)) = (0,\ y + c'_1\left(\beta \vee \beta^{-1}\right)+ \beta)$, then from the evolution equation of the diffusion in~\eqref{eq:diffusionjsq}, the event $\cE_t^2$ implies the event
$$\tilde{\cE}_t^2:=\Bigg[Q_1(t) + Q_2(t) \leq y+ c'_1\left(\beta \vee \beta^{-1}\right) + \beta +\sqrt{2}W(t) - \frac{\beta t}{2},  \inf_{s\leq t} Q_2(s) > c'_1\left(\beta \vee \beta^{-1}\right)+ \beta\Bigg].$$
Therefore,
\begin{align}\label{eq:lem3.7-1}
\prob_{(0,\ y + c'_1\left(\beta \vee \beta^{-1}\right)+ \beta)}\big(\cE_t\big) 
&\leq \prob_{(0,\ y + c'_1\left(\beta \vee \beta^{-1}\right)+ \beta)}\Big(\cE_t^1\Big)+\prob_{(0,\ y + c'_1\left(\beta \vee \beta^{-1}\right)+ \beta)}\Big(\tilde{\cE}_t^2\Big).
\end{align}
Now, choose $c'_1, c'_2$ as in Lemma \ref{lem:q1integral}. Then for any $y \ge 1$,
\begin{equation}\label{eq:lem3.7-2}
\prob_{(0,\ y + c'_1\left(\beta \vee \beta^{-1}\right)+ \beta)}\Big(\cE_t^1\Big) \leq \exp(-c'_2t^{1/5}\left(\beta \vee \beta^{-1}\right)^{2/5}).
\end{equation}
Also, note that 
\begin{equation}\label{eq:lem3.7-3}
\begin{split}
&\prob_{(0,\ y + c'_1\left(\beta \vee \beta^{-1}\right)+ \beta)}\Big(\tilde{\cE}_t^2\Big)\\
&\leq \prob_{(0,\ y + c'_1\left(\beta \vee \beta^{-1}\right)+ \beta)}\big(Q_1(t) \leq y+\sqrt{2}W(t) - \frac{\beta t}{2},  \inf_{s\leq t} Q_2(s) > c'_1\left(\beta \vee \beta^{-1}\right)+ \beta\big)\\
&\leq \Pro{\sqrt{2}W(t)>\frac{\beta t}{4}} + \prob_{(0,\ y + c'_1\left(\beta \vee \beta^{-1}\right)+ \beta)}\left(Q_1(t) \leq y - \frac{\beta t}{4},  \inf_{s\leq t} Q_2(s) > c'_1\left(\beta \vee \beta^{-1}\right)+ \beta\right).
\end{split}
\end{equation}
Due to Brownian scaling we have
\begin{equation}
\Pro{\sqrt{2}W(t)>\frac{\beta t}{4}}\leq c\exp(-c'\beta^2 t) \ \text{ for } \ t \ge \beta^{-2},
\end{equation}
where $c,c'$ do not depend on $\beta$.
Moreover, choosing $t > 8y/\beta$, and applying Lemma~\ref{lem:integrated}~\eqref{lem:tails} and Lemma~\ref{lem:integrated}~\eqref{lem:Nt-upper} with $\varepsilon = (\beta \wedge \beta^{-1})/4$ and $M=c'_1\left(\beta \vee \beta^{-1}\right)$,
\begin{equation}\label{eq:lem3.7-last}
\begin{split}
&\prob_{(0,\ y + c'_1\left(\beta \vee \beta^{-1}\right)+ \beta)}\big(Q_1(t) \leq y - \frac{\beta t}{4},  \inf_{s\leq t} Q_2(s) > c'_1\left(\beta \vee \beta^{-1}\right)+ \beta\big)\\
&\leq \prob_{(0,\ y + c'_1\left(\beta \vee \beta^{-1}\right)+ \beta)}\big(Q_1(t) \leq - \frac{\beta t}{8},  \inf_{s\leq t} Q_2(s) > c'_1\left(\beta \vee \beta^{-1}\right)+ \beta\big)\\
&\leq \prob\Big(\sup_{1\leq i\leq N_t}u_i > \frac{\beta t}{8}\Big)\\
&\leq \Pro{N_t>16b\left(\beta \vee \beta^{-1}\right)^2t} + 16b\left(\beta \vee \beta^{-1}\right)^2t \Pro{u_1>\frac{\beta t}{8}}\\
&\leq \exp(-c\left(\beta \vee \beta^{-1}\right)^2t) + 16b\left(\beta \vee \beta^{-1}\right)^2t\ \exp\left(-\left(\beta \vee \beta^{-1}\right)\left(\frac{\beta t}{8}-\frac{\beta \wedge \beta^{-1}}{4})\right)\right)\\
&\leq \exp(-c\left(\beta \vee \beta^{-1}\right)^2t) + 16b\left(\beta \vee \beta^{-1}\right)^2t\ \exp\left(-\left(\beta \vee \beta^{-1}\right)\left(\frac{\beta t}{16}\right)\right),
\end{split}
\end{equation}
where $b,c$ do not depend on $\beta$. 
Combining Equations~\eqref{eq:lem3.7-1} -- \eqref{eq:lem3.7-last} completes the proof of the lemma.
\end{proof}

\section{Proof of Lemma~\ref{lem:alpha}}\label{app:second}
In order to prove Lemma~\ref{lem:alpha}, set $M>0$ to be a fixed large number to be chosen later and $(Q_1(0), Q_2(0))=(x,y)$ for some $x \in [-M/2,0], y \in (0,B]$.
For $i\geq 1$ define the stopping times
\begin{align*}
\tau_{2,2i-1}&:= \inf\Big\{t\geq \tau_{2,2i-2}: Q_2(t)= 2B\quad\mbox{or}\quad Q_1(t)=-M\Big\},\\ 
\tau_{2,2i}&:= \inf\Big\{t\geq \tau_{2,2i-1}: Q_2(t)= 2B \quad\mbox{or}\quad Q_1(t)=-\frac{M}{2}\Big\},
\end{align*}
where by convention we take $\tau_{2,0}\equiv 0$.
Also define 
$$N^*:= \inf \Big\{k\geq 0: Q_2(\tau_{2,2k+1})= 2B\Big\}.$$
Therefore, note that 
\begin{equation}\label{eq:tau2B}
\tau_2(2B) = \sum_{j=1}^{2N^*+1}(\tau_{2,j}-\tau_{2,j-1}).
\end{equation}
The proof of Lemma \ref{lem:alpha} consists of three parts:
(i) Lemma~\ref{lem:q2upcross} contains the required probability estimate to analyze the time interval $\tau_{2, 2i-1}- \tau_{2,2i-2}$,
(ii) Lemma~\ref{lem:evencycle} estimates the tail probabilities for the time interval $\tau_{2, 2i}- \tau_{2,2i-1}$,
and  (iii) Lemma~\ref{lem:N-star-tail} provides the tail probabilities for the random variable $N^*$.
Lemma~\ref{lem:exit-by-q2} is used in the proof of Lemma~\ref{lem:N-star-tail}.
Combining Equation~\eqref{eq:tau2B} and Lemmas \ref{lem:q2upcross}, \ref{lem:evencycle}, and \ref{lem:N-star-tail}, we will complete the proof of Lemma~\ref{lem:alpha}.

\begin{lemma}\label{lem:q2upcross}
For any fixed $B,M> 0$, there exists $p^{(1)}(M,B)>0$, such that
$$\inf_{\substack{x\in[-M,\ 0],\\ y \in (0, 2B]}}\prob_{(x,y)}\Big(\sup_{0\leq s\leq 1}Q_2(s)>2B\Big)\geq p^{(1)}(M,B).$$
\end{lemma}
\begin{proof}
Recall that 
$$Q_1(t) = Q_1(0) +\sqrt{2}W(t) -\beta t +\int_0^t (-Q_1(s) + Q_2(s))\dif s - L(t),$$
where 
\begin{equation}\label{eq:lt}
\begin{split}
L(t) &=\ \sup_{s\leq t} \Big(Q_1(0)+\sqrt{2}W(s)-\beta s + \int_0^s(-Q_1(u)+Q_2(u))\dif u\Big)^+
\geq\ \sup_{s\leq t} (Q_1(0)+\sqrt{2}W(s)-\beta s)^+.
\end{split}
\end{equation}
Thus, $\Pro{L(1)> 4B}\geq \Pro{\sqrt{2}W(1)>\beta + 4B -Q_1(0)}$.
Observe that for any $Q_2(0) = y \le 2B$,
$$\big\{L(1)>4B\big\}\implies \big\{sup_{s\leq 1}Q_2(s)>2B\big\}$$
To see this, suppose $L(1)>4B$. If $\sup_{s\leq 1} Q_2(s) \leq 2B$, then
$$Q_2(1) = y + L(1) - \int_0^1 Q_2(s) \dif s \ge L(1) - 2B > 2B$$
which is a contradiction.
Therefore, 
\begin{align*}
&\inf_{\substack{x\in[-M,\ 0],\\ y \in (0,2B]}}\prob_{(x,y)}\Big(\sup_{0\leq s\leq 1}Q_2(s)>2B\Big)
\geq
\inf_{\substack{x\in[-M,\ 0],\\ y \in (0, 2B]}}\prob_{(x,y)}\Big(L(1)>4B\Big)\\
&\geq
\inf_{\substack{x\in[-M,\ 0],\\ y \in (0, 2B]}}\prob_{(x,y)}\Big(\sqrt{2}W(1)>\beta + 4B - x\Big)\\
&\geq \Pro{\sqrt{2}W(1)>\beta + 4B+M} = p^{(1)}(M,B)>0.
\end{align*}
This completes the proof of Lemma~\ref{lem:q2upcross}.
\end{proof}

\begin{lemma}\label{lem:evencycle}
For any $j \ge 0$ and any fixed $M\ge 6\beta$, there exists $c_\tau^{(1)}>0$ such that for all $t\geq 2$,
$$\sup_{\substack{x\in[-M/2,\ 0],\\ y \in (0, B]}}\prob_{(x,y)}\Big(\tau_{2,2j+2} - \tau_{2, 2j+1}>t\ \Big|\ N^*>j\Big)\leq \exp(-c_\tau^{(1)} t).$$
\end{lemma}
\begin{proof}
Let us denote $Q_1^* = Q_1+\beta$. Since $ N^*>j$, we know $Q_2(\tau_{2, 2j+1}) < 2B$.
In that case, for $t>\tau_{2,2j+1},$
\begin{align*}
Q_1^*(t) &= Q_1^*(\tau_{2, 2j+1}) + \sqrt{2}W(t) + \int_{\tau_{2,2j+1}}^t (-Q_1^*(s)+Q_2(s))\dif s\\
&\geq Q_1^*(\tau_{2, 2j+1}) +\sqrt{2}W(t) - \int_{\tau_{2,2j+1}}^tQ_1^*(s)\dif s\\
&= -M+\beta + \sqrt{2}W(t) - \int_{\tau_{2,2j+1}}^t Q_1^*(s)\dif s.
\end{align*}
Thus, we obtain
\begin{align*}
\prob_{(x,y)}\Big(\tau_{2,2j+2} - \tau_{2, 2j+1}>t\ \Big|\ N^*>j\Big)
\leq \prob\Big(\sup_{s\leq t}(\sqrt{2}W(s) - (-M/2+\beta)s)\leq M/2\Big),
\end{align*}
since for $t\in (\tau_{2,2j+1}, \tau_{2,2j+2})$, $Q_1^*(s)\leq -M/2+\beta$.
Therefore, as $M \ge 6\beta$, for all $t\geq 2$,
\begin{align*}
\prob_{(x,y)}\Big(\tau_{2,2j+2}-\tau_{2,2j+1}>t\ \Big|\ N^*>j\Big)
&\leq \Pro{\sqrt{2}W(t) \leq M/2- (M/2-\beta)t}\\
&\le \Pro{\sqrt{2}W(t) \leq - (M/2-\beta)t/4}\\
&\leq \exp(-c_\tau^{(1)}(M/2-\beta)^2 t)
 \le \exp(-c_\tau^{(1)} t),
\end{align*}
where $c_\tau^{(1)}$ does not depend on $x,y$.
\end{proof}

\begin{lemma}\label{lem:exit-by-q2}
For any fixed $B>0$ and $M>8B + 2\beta$, there exists $p^{(2)}=p^{(2)}(M,B)>0$ such that
\begin{align*}
\inf_{\substack{ x\in [-M/2,\ 0], \\ y \in (0, B]}}\prob_{(x,y)}\Big(\exists\ t^*\in [0,1],\ \mbox{\rm such that } \sup_{0\leq t\leq t^*}Q_2(t) \ge 2B, \inf_{0\leq t\leq t^*}Q_1(t)>-M\Big)
\geq p^{(2)}.
\end{align*}
\end{lemma}
\begin{proof}
For fixed $B>0$ and $M>8B + 2\beta$, consider the event 
$$\cE(\beta,M):= \Big\{\sqrt{2}W(1)>\beta + 4B + \frac{M}{2}, \quad \inf_{t\in [0,1]}\sqrt{2}W(t)>\beta + 4B-\frac{M}{2}\Big\}.$$
From the representation \eqref{eq:lt}, note that the event $\cE(M,B)$ implies the event 
$\{L(1) > 4B\}$, which in turn implies that there exists $t^*\in [0,1]$ such that  $L(t^*)=4B$
and $\forall\ t\leq t^*$,
\begin{align*}
Q_1(t) &\geq -\frac{M}{2} +\sqrt{2}W(t) -\beta -4B
>-\frac{M}{2} -\beta -4B + \big(\beta + 4B -\frac{M}{2}\big)=-M.
\end{align*}
Therefore, $\inf_{0\leq t\leq t^*}Q_1(t) > -M$.
Furthermore, we claim that $\sup_{0\leq t\leq t^*}Q_2(t) \ge 2B$.
Indeed, if $\sup_{0\leq t\leq t^*}Q_2(t) < 2B$, then
$$
Q_2(t^*)\geq L(t^*) - \int_0^{t^*}Q_2(s)\dif s >4B - 2B t^*\geq 2B,
$$
since $0\leq t^*\leq 1$, which leads to a contradiction.
Finally, 
\begin{multline*}
\inf_{\substack{ x\in [-M/2,\ 0], \\ y \in (0, \beta^{-1}]}}\prob_{(x,y)}\Big(\exists\ t^*\in [0,1],\ \mbox{\rm such that } \sup_{0\leq t\leq t^*}Q_2(t)>2B, \inf_{0\leq t\leq t^*}Q_1(t)>-M\Big)\\
\geq \Pro{\cE(M,B)} >0.
\end{multline*}
This completes the proof of the lemma.
\end{proof}

\begin{lemma}\label{lem:N-star-tail}
For any fixed $B>0$ and $M> 8B + 2\beta$, there exist $c_N^{(2)}, n_N>0$ such that for all $n\geq n_N$, 
$$\sup_{\substack{x\in[-M/2,\ 0],\\ y \in (0, B]}}\prob_{(x,y)}(N^*>n)\leq \exp(-c_N^{(2)}n).$$
\end{lemma}
\begin{proof}
Observe that
\begin{align*}
\prob_{(x,y)}(N^*>n)\leq \prob_{(x,y)}(Q_1(\tau_{2,2k+1})=-M \mbox{ and }Q_2(\tau_{2,2k+1})< 2B \text{ for all } k \le n)\leq (1-p^*)^n,
\end{align*}
using strong Markov property, where 
\begin{align*}
p^* &:= \inf_{\substack{x\in[-M/2,\ 0]\\ y\in(0, B]}}\prob_{(x,y)}(Q_2\mbox{ hits } 2B \mbox{ before }Q_1\mbox{ hits }-M)\\
&\geq \inf_{\substack{x\in[-M/2,\ 0]\\ y\in(0, B]}}\prob_{(x,y)}(\exists\ t^*\in[0,1]\mbox{ such that }\sup_{0\leq t\leq t^*}Q_2(t) >2B, \inf_{0\leq t\leq t^*}Q_1(t)>-M)\\
&\geq p^{(2)}(M,B)>0,
\end{align*}
by Lemma \ref{lem:exit-by-q2}, choosing $M> 8B + 2\beta$.
\end{proof}

Now, we have all the necessary results to prove Lemma \ref{lem:alpha}.
\begin{proof}[Proof of Lemma~\ref{lem:alpha}]
Recall that $\tau_2(2B) = \sum_{j=1}^{2N^*+1}(\tau_{2,j}-\tau_{2,j-1}).$
From Lemma~\ref{lem:q2upcross} observe that for any fixed $M> 0$ and any $x \in [-M/2,0], \ y \in (0,B]$,
\begin{align*}
\prob_{(x, y)}(\tau_{2,1}>n) &=\mathbb{E}_{(x, y)}\left(\mathbbm{1}_{[\tau_{2,1}>n-1]}\prob_{(Q_1(n-1), Q_2(n-1))}(\tau_{2,1}>1)\right)\\
&\leq (1-p^{(1)}(M,B))\Pro{\tau_{2,1}>n-1},
\end{align*}
Which implies $\prob_{(x,y)}(\tau_{2,1}>n) \leq (1-p^{(1)}(M,B))^n$.
Furthermore, following the same argument as above, we can claim that for all $j\geq 1$,
\begin{equation}\label{eq:oddcycle}
\prob_{(x,y)}(\tau_{2,2j-1}-\tau_{2,2j-2}\geq n)\leq (1-p^{(1)}(M,B))^n.
\end{equation}
Therefore for $t\geq 9$, choosing $M> 8B + 6\beta$, we can write for any $x \in [-M/2,0], \ y \in (0,B]$,
\begin{align*}
&\prob_{(x,y)}(\tau_2(2B)>t)
\leq \prob_{(x,y)}(N^*>n) + \prob_{(x,y)}\Big(\sum_{j=1}^{2n+1}(\tau_{2,j}-\tau_{2,j-1})>t\Big) \\
&\leq \exp(-c_N^{(2)}n) + (2n+1)\exp(-ct/(2n+1)),\hspace{1cm} \mbox{Due to Lemmas~\ref{lem:evencycle} and~\ref{lem:N-star-tail}, and~\eqref{eq:oddcycle}}\\
&\leq c'\sqrt{t}\e^{-c\sqrt{t}}
\leq \e^{c_\alpha^{(2)}\sqrt{t}},\hspace{4.5cm}\mbox{[choosing }n = \lfloor(\sqrt{t}-1)/2\rfloor\mbox{]}
\end{align*} 
where $c_\alpha^{(2)}$ does not depend on $(x,y)$.
\end{proof}

\end{appendices}

\end{document}